\documentclass[12pt,reqno,twoside,a4paper]{amsart}

\usepackage{amsmath, amssymb, amscd, amsthm}

\theoremstyle{plain}
\newtheorem{thm}{Theorem}[section]
\newtheorem{cor}[thm]{Corollary}

\newtheorem{lem}[thm]{Lemma}
\newtheorem{prop}[thm]{Proposition}
\newtheorem{Def}[thm]{Definition}

\theoremstyle{remark}
\newtheorem{rem}[thm]{Remark}
\numberwithin{equation}{section}

\newcommand{\R}{\mathbb R}

\newcommand{\C}{\mathbb C}
\newcommand{\Z}{\mathbb Z}
\newcommand{\T}{\mathbb T}
\newcommand{\al}{\alpha}
\newcommand{\be}{\beta}
\newcommand{\ga}{\gamma}
\newcommand{\Ga}{\Gamma}
\newcommand{\de}{\delta}

\newcommand{\eps}{\varepsilon}
\newcommand{\si}{\sigma}
\newcommand{\te}{\theta}

\newcommand{\la}{\lambda}

\newcommand{\hf}{\frac{1}{2}}
\newcommand{\ph}[5]{\,_{#1}\phi_{#2} \left( \genfrac{.}{.}{0pt}{}{#3}{#4}\ ;#5 \right)}
\newcommand{\Res}[1]{\underset{#1}{\mathrm{Res}}\,}
\newcommand{\vect}[2]{\left( \genfrac{.}{.}{0pt}{}{#1}{#2} \right) }

\begin{document}
\title{The vector-valued big $q$-Jacobi transform}

\author{Wolter Groenevelt}
\address{Korteweg-De Vries Instituut voor Wiskunde\\
Universiteit van Amsterdam\\
Plantage Muidergracht 24\\
1018 TV Amsterdam\\
The Netherlands}
\email{wgroenev@science.uva.nl}
\thanks{The author is supported by the Netherlands Organization for Scientific Research (NWO) for the Vidi-project ``Symmetry and modularity in exactly solvable models.''}

\begin{abstract}
Big $q$-Jacobi functions are eigenfunctions of a second order $q$-difference operator $L$.
We study $L$ as an unbounded self-adjoint operator on an $L^2$-space of functions on $\mathbb R$ with a discrete measure. We describe explicitly the spectral decomposition of $L$ using an integral transform $\mathcal F$ with two different big $q$-Jacobi functions as a kernel, and we construct the inverse of $\mathcal F$.
\end{abstract}

\maketitle

\section{Introduction}
Integral transforms with a hypergeometric function as a kernel have been the subject of many papers in the literature. A famous example is the Jacobi transform, first studied by Weyl \cite{We10}, which is an integral transform with a certain $_2F_1$-function, the Jacobi function, as a kernel. The inverse of the Jacobi transform can be obtained from spectral analysis of the hypergeometric differential operator $D$, which is an unbounded self-adjoint operator on a weighted $L^2$-space of functions on $[0,\infty)$. We refer to \cite{Koo84} for a survey on Jacobi functions. In a recent paper \cite{Ne05} Neretin studied the hypergeometric differential operator $D$ as a self-adjoint operator on a weighted $L^2$-space of functions on $\R$. In this setting the spectral analysis of $D$ leads to an integral transform with two different Jacobi functions (vector-valued Jacobi functions) as a kernel, corresponding to the multiplicity two of the continuous spectrum.   

In this paper we obtain a $q$-analogue of Neretin's vector-valued Jacobi transform (or double index hypergeometric transform). There exist several $q$-analogues of Jacobi functions, namely the little and big $q$-Jacobi functions and the Askey-Wilson functions, see \cite{K95}, \cite{KMU95}, \cite{KS01}, \cite{KS03}. Here we consider the big $q$-Jacobi function, which is a basic hypergeometric $_3\phi_2$-function that is the kernel in the big $q$-Jacobi transform by Koelink and Stokman \cite{KS03}. The big $q$-Jacobi transform and its inverse arise from spectral analysis of a second order $q$-difference operator $L$, that is an unbounded self-adjoint operator on an $L^2$-space consisting of square integrable functions with respect to a discrete measure on $[-1,\infty)$. In this paper we study the same $q$-difference operator $L$ as an unbounded self-adjoint operator on a different Hilbert space, namely an $L^2$-space of functions on $\R$ with a discrete measure. The continuous spectrum of $L$ has multiplicity two, thus leading to an integral transform pair with two different big $q$-Jacobi functions as a kernel.  We call this the vector-valued big $q$-Jacobi transform.

The vector-valued Jacobi transform has an interpretation in the representation theory of Lie algebra $\mathfrak{su}(1,1)$ (or equivalently $\mathfrak{sl}(2,\R)$) as follows, see \cite[Section 4]{Ne05}. The hypergeometric differential operator $D$ arises from a suitable restriction of the Casimir operator in the tensor product of two principal unitary series. The spectral analysis of $D$ now gives the decomposition into irreducible representations, and the vector-valued Jacobi transform can be used to construct explicitly the intertwiner for these representations. The multiplicity two of the continuous spectrum corresponds to the multiplicity of the principal unitary series occurring in the decomposition, see \cite{Pu61}, \cite{Ma75} for the precise decomposition. There is a similar interpretation of the vector-valued big $q$-Jacobi transform in the representation theory of the quantized universal enveloping algebra $\mathcal U_q(\mathfrak{su}(1,1))$. However, the corresponding representation is no longer a tensor product representation, but a sum of two tensor products of principal unitary series. This will be the subject of a future paper.

The big $q$-Jacobi functions are nonpolynomial extensions of the big $q$-Jacobi polynomials \cite{AA85}, but they can also be considered as extensions of the continuous dual $q^{-1}$-Hahn polynomials, see \cite{KS03}. In this light, the vector-valued big $q$-Jacobi transform may also be considered as a $q$-analogue of the integral transform corresponding to the $_3F_2$-functions $(\Xi_n^{(1)}, \Xi_n^{(2)})$ from \cite[Thm.1.3]{Ne05}, and of the continuous Hahn transform from \cite{GKR05}. In both transforms the kernel consists of two $_3F_2$-functions that are extensions of the continuous dual Hahn polynomials.

The organization of this paper is as follows. In Section \ref{sec2} we introduce the second order $q$-difference operator $L$ and a weighted $L^2$-space of functions on $\R_q$, a $q$-analogue of the real line. The difference operator $L$ is an unbounded operator on this $L^2$-space. We define the Casorati determinant, a difference analogue of the Wronskian, and with the Casorati determinant we determine a dense domain on which $L$ is self-adjoint. In section \ref{sec3} we introduce the big $q$-Jacobi functions as eigenfunctions of $L$ given by a specific $_3\phi_2$-series. We also give the asymptotic solutions, which are $_3\phi_2$-series with nice asymptotic behavior at $+\infty$ or $-\infty$. A crucial point here is the fact that all eigenfunctions that we consider can uniquely be extended to functions on $\R_q$. In Section \ref{sec4} we define the Green kernel using the asymptotic solutions and we determine the spectral decomposition for $L$. In Section 5 we define the vector-valued big $q$-Jacobi transform $\mathcal F$, and we determine its inverse. A left inverse $\mathcal G$ of $\mathcal F$ follows immediately from the spectral analysis done in Section \ref{sec4}. To show that $\mathcal G$ is also a right inverse, we use a classical method that essentially comes down to approximating with the Fourier transform. Finally, in the appendix two lemmas are proved which involve rather long computations. \\

\textit{Notations}. 
We use standard notations for $q$-shifted factorials, $\te$-functions and basic hypergeometric series \cite{GR04}. 
We fix a number $q \in (0,1)$. The $q$-shifted factorials are defined by
\[
(x;q)_\infty = \prod_{k=0}^\infty (1-xq^k), \qquad (x;q)_n = \frac{ (x;q)_\infty }{(xq^n;q)_\infty}, \ n \in \Z.
\]
The (normalized) Jacobi $\te$-function is defined by
\[
\te(x) = (x,q/x;q)_\infty, \qquad x \not\in q^\Z.
\]
From this definition it follows that the $\te$-function satisfies
\[
\te(x)=\te(q/x) =-x\te(qx) = -x \te(1/x). 
\]
We often use these identities without mentioning them. For products of $q$-shifted factorials and products of $\te$-functions we use the shorthand notations 
\begin{align*}
(x_1,x_2,\ldots, x_k;q)_n &=(x_1;q)_n (x_2;q)_n \cdots (x_k;q)_n , & n \in \Z \cup \{\infty \},\\
\te(x_1,x_2,\ldots,x_k)&=\te(x_1)\te(x_2)\cdots \te(x_k),
\end{align*}
and
\[
(xy^{\pm 1};q)_\infty= (xy,x/y;q)_\infty, \qquad \te(xy^{\pm 1})= \te(xy,x/y).
\]
An identity for $\te$-functions that we frequently use is
\begin{equation} \label{eq:theta prod}
\te(xv,x/v,yw,y/w)-\te(xw,x/w,yv,y/v) = \frac{y}{v} \te(xy,x/y,vw,v/w),
\end{equation}
see \cite[Exer. 2.16(i)]{GR04}. The basic hypergeometric function $_r\phi_s$ is defined by
\[
\ph{r}{s}{x_1,x_2,\ldots, x_r}{y_1,y_2,\ldots,y_s}{q,z} = \sum_{k=0}^\infty \frac{ (x_1,x_2,\ldots,x_r;q)_k }{ (q,y_1,y_2,\ldots,y_s;q)_k } \Big((-1)^k q^{k(k-1)/2} \Big)^{1+s-r} z^k.
\]

\section{The second order $q$-difference operator} \label{sec2}
In this section we introduce an unbounded second order $q$-difference operator $L$ acting on functions on a $q$-analogue of the real line, and we determine a dense domain on which $L$ is self-adjoint.
\subsection{The difference operator}We fix two real numbers $z_+>0$ and $z_-<0$. Let $\R_q^+$ and $\R_q^-$ be the two sets
\[
\R_q^+= \{z_+q^n\ | \ n \in \Z \}, \qquad \R_q^-= \{z_-q^n\ | \ n \in \Z \},
\]
and define 
\[
\R_q=\R_q^- \cup \R_q^+,
\]
which we consider as a $q$-analogue of the real line. For $x \in \R_q$ we sometimes write $x=zq^{k}$, which means that $z=z_-$ or $z=z_+$, and $k \in \Z$. We denote by $F(\R_q)$ the linear space of complex-valued functions on $\R_q$. 

The second order difference operator $L$ we are going to study depends on four parameters. Let $P_{q,z_-,z_+}$ be the set consisting of pairs of parameters $(\al,\be) \in \C^2$ such that $\al,\be \not\in z_{\pm}^{-1}q^\Z$, and one of the following conditions is satisfied:
\begin{itemize}
\item $\al=\overline{\be}$
\item there exists a $k_0 \in \Z$ such that $z_+q^{k_0} <\be^{-1}<\al^{-1}< z_+q^{k_0-1}$
\item there exists a $k_0 \in \Z$ such that $z_-q^{k_0-1} <\al^{-1}<\be^{-1}< z_-q^{k_0}$
\end{itemize}
In particular this implies that $q <|\al/\be|\leq 1$, and $\al$ and $\be$ have the same sign in case they are real.
We define the parameter domain $P$ to be the following set,
\[
P = \big\{ (a,b,c,d) \in \C^4 \ | \ (a,b) \in P_{q,z_-,z_+},\, (c,d) \in P_{q,z_-,z_+},\, a\neq b \big\}.
\]
From here on we assume that $(a,b,c,d) \in P$, unless explicitly stated otherwise.

We define a linear operator $L=L_{a,b,c,d}: F(\R_q) \rightarrow  F(\R_q)$ by
\[
L = A(\cdot) T_{q^{-1}}  + B(\cdot) T_q +C(\cdot)id,
\]
where $T_\al$ is the shift operator $(T_\al f)(x)= f(\al x)$ for $\al \in \C$, $id$ denotes the identity operator, and
\[
\begin{split}
A(x) &= s^{-1}\left(1-\frac{q}{ax}\right) \left(1-\frac{q}{bx}\right),\\ 
B(x) &= s\left(1-\frac{1}{cx}\right)\left(1-\frac{1}{dx}\right), \\
C(x)&= s^{-1}+s-A(x)-B(x),
\end{split}
\]
where $s=\sqrt{cdq/ab}$. Here we use the usual branch of $\sqrt{\cdot}$ that is positive on $\R_+$. Note that the conditions on the parameters ensure that $A(x) \neq 0$ and $B(x) \neq 0$ for all $x \in \R_q$.
\begin{rem} \label{rem:symmetry}
(a) There is an obvious symmetry in the parameters;
\[
L_{a,b,c,d} = L_{b,a,c,d} = L_{a,b,d,c}.
\]
This will be useful when we study eigenfunctions of $L$ later on.

(b) Let $a=z_+^{-1}q^{1-m}$, $m \in \Z$ fixed, (so $(a,b,c,d) \not\in P$), then the coefficient $A(x)$ vanishes at the point $x=z_+q^m$. In this case certain restrictions of the difference operator $L$ are well known in the literature. Let $L_-$ (respectively $L_+$) denote the operator $L$ restricted to functions on $\{z_+q^k\ | \ k \in \Z_{\geq m} \}$ (respectively $\{z_+q^k\ | \ k \in \Z_{\leq m-1} \}$). Then $L_\pm$ is equivalent to the Jacobi operator for the continuous dual $q^{\pm 1}$-Hahn polynomials, which are Askey-Wilson polynomials with one of the parameters equal to zero,  see \cite{AW85}, \cite{KS98}. Moreover, the operator $L$ restricted to functions on $\R_q^- \cup \{z_+q^k\ | \ k \in \Z_{\geq m} \}$ is equivalent to the difference operator studied by Koelink and Stokman \cite{KS03} to obtain the big $q$-Jacobi function transform.

If we also set $c=z_-^{-1}q^{-n}$, $n \in \Z$ fixed, then the coefficient $B(x)$ vanishes at the point $x=z_-q^n$. In this case the operator $L$ restricted to functions on $\{z_-q^k \ | \ k \in \Z_{\geq n} \} \cup \{z_+q^k\ | \ k \in \Z_{\geq m} \}$ is equivalent to the difference operator for big $q$-Jacobi polynomials, see \cite[Sect.10]{KS03} where the orthogonality relations for the big $q$-Jacobi polynomials are obtained from spectral analysis of $L$. See also \cite{AA85} and \cite{KS98} for the big $q$-Jacobi polynomials.
\end{rem}

\subsection{The Casorati determinant}
The Jackson $q$-integral is defined by
\begin{gather*}
\int_{0}^\al f(x) d_qx =  (1-q) \sum_{k=0}^\infty f(\al q^k)\al q^k,\\
\int_\al^\be f(x) d_qx =  \int_0^\be f(x) d_qx - \int_0^\al f(x) d_qx\\
\int_0^{\infty(\al)} f(x) d_qx = (1-q) \sum_{k=-\infty}^\infty f(\al q^k) \al q^k,
\end{gather*}
for $\al,\be \in \C^*$, and $f$ is a function such that the sums converge absolutely. We will denote
\[
\int_{\R_q} f(x) d_qx = \int_0^{\infty(z_+)} f(x) d_qx - \int_0^{\infty(z_-)} f(x) d_qx,
\]
We define a weight function $w$ on $\R_q$ by
\begin{equation} \label{def:weight w}
w(x)=w(x;a,b,c,d;q) = \frac{ (ax, bx;q)_\infty }{ (cx, dx;q)_\infty } .
\end{equation}
Note that for $(a,b,c,d) \in P$ the weight function $w$ is positive on $\R_q$, and $w$ is continuous at the origin. Let $\mathcal L^2=\mathcal L^2(\R_q,w(x)d_qx)$ be the Hilbert space consisting of functions $f \in F(\R_q)$ that have finite norm with respect to the inner product
\[
\langle f, g \rangle_{\mathcal L^2}=\int_{\R_q} f(x) \overline{g(x)} w(x) d_qx.
\]
For $k,l,m,n \in \Z$ we define a truncated inner product by
\[
\langle f, g \rangle_{k,l;m,n} =
 \int_{z_-q^k}^{z_-q^{l+1}} f(x) \overline{g(x)}w(x) d_qx + \int_{z_+q^{m+1}}^{z_+q^n} f(x) \overline{g(x)}w(x) d_qx \\
\]
If $f,g \in \mathcal L^2$ we have
\[
\lim_{\substack{l,m \rightarrow \infty \\k,n \rightarrow -\infty}} \langle f, g \rangle_{k,l;m,n} = \langle f, g \rangle_{\mathcal L^2}.
\]
We define a function $u$ on $\R_q$, closely related to the weight function $w$, by
\[
u(x) = (1-q)^2B(x)x^2w(x) = (1-q)^2 \sqrt{q/abcd}  \frac{ (ax,bx;q)_\infty }{ (cqx,dqx;q)_\infty }.  
\]\

\begin{Def} \label{def:Casorati}
For $f,g \in F(\R_q)$ we define the Casorati determinant $D(f,g) \in F(\R_q)$ by
\[
\begin{split}
D(f,g)(x)&= \Big( f(x) g(qx) - f(qx) g(x) \Big) \frac{u(x)}{(1-q)x} \\
&= \Big( (D_qf)(x) g(x) - f(x) (D_qg)(x) \Big) u(x).
\end{split}
\]
\end{Def}
Here $D_q:F(\R_q)\rightarrow F(\R_q)$ is the $q$-difference operator given by
\[
(D_qf)(x) = \frac{ f(x) - f(qx) }{x(1-q)}.
\]

\begin{prop} \label{prop:Casorati}
For $f,g \in F(\R_q)$ we have
\[
\begin{split}
\langle Lf, g \rangle_{k,l;m,n} &- \langle f,L g \rangle_{k,l;m,n} =\\
& D(f,\overline{g})(z_-q^{l})-D(f,\overline{g})(z_-q^{k-1})+ D(f,\overline{g})(z_+q^{n-1}) -D(f,\overline{g})(z_+q^{m}).
\end{split}
\]
\end{prop}
\begin{proof}
For $f,g \in F(\R_q)$,
\[
\begin{split}
\Big(&(Lf)(x)g(x) - f(x) (Lg)(x)\Big) xw(x)  =\\
 &A(x) \Big( f(x/q)g(x) - f(x)g(x/q) \Big) x w(x) \\
& \quad+ B(x)\Big( f(xq)g(x) - f(x)g(xq) \Big) xw(x) \\
=& \Big( f(x/q)g(x) - f(x)g(x/q) \Big)  \frac{ (ax/q,bx/q;q)_\infty } { (cx,dx;q)_\infty }  \frac{q^2}{xsab}  \\
& \quad + \Big( f(xq)g(x) - f(x)g(xq) \Big) \frac{ (ax,bx;q)_\infty }{ (cqx,dqx;q)_\infty } \frac{s}{cdx} 
\end{split}
\]
Note that $abs/q =\sqrt{abcd/q}= cd/s$, so we obtain
\[
\Big((Lf)(x)g(x) - f(x) (Lg)(x)\Big) (1-q)xw(x) = D(f,g)(x/q)- D(f,g)(x).
\]
Now the sums of the truncated inner products in the lemma become telescoping, and then the result follows.
\end{proof}
In order to determine a suitable domain on which $L$ is self-adjoint, we need to find the limit behavior of Casorati determinants. First, to find the asymptotic behavior of $D(f,g)(x)$ for large $x$, we need the behavior of $u(x)$ for large $x$.
\begin{lem} \label{lem:asymp w&u}
Let $x=zq^{-k} \in \R_q$, then for $k \rightarrow \infty$
\[
\begin{split}
xw(x) &=  (1-q)^{-1}K_zs^{-2k} \big(1+\mathcal O(q^k) \big), \\
\frac{u(x)}{x} &= (1-q)K_z  s^{1-2k}\big(1+\mathcal O(q^k) \big),
\end{split}
\]
where 
\[
K_z =K_z(a,b,c,d;q)= z(1-q) \frac{ \te(az,bz)}{ \te(cz,dz)}.
\]

\end{lem}
\begin{proof}
Let $x=zq^{-k} \in \R_q$. Using the identity
\[
\frac{ (\al q^{-n};q)_n }{ (\be q^{-n};q)_n} = \frac{ (q/\al;q)_n }{ (q/\be;q)_n } \left( \frac{\al}{\be} \right)^n,
\]
and the definition \eqref{def:weight w} of the weight function $w$ we obtain 
\[
w(zq^{-k}) =  \frac{ (q/az,q/bz;q)_k (az,bz;q)_\infty }{ (q/cz,q/dz;q)_k (cz,dz;q)_\infty } \left(\frac{ab}{cd}\right)^{k},
\]
and
\[
u(zq^{-k}) = (1-q)^2 \sqrt{q/abcd} \frac{ (q/az,q/bz;q)_k (az,bz;q)_\infty} { (1/cz,1/dz;q)_k (zcq,zdq;q)_\infty }   \left(\frac{ab}{cdq^2}\right)^{k}
\]
From this the asymptotic behavior of $xw(x)$ and $u(x)/x$ for large $x$ follows.
\end{proof}

\begin{lem} \label{lem:lim infty}
Let $f,g \in \mathcal L^2$, then
\[
\lim_{x \rightarrow \pm \infty} D(f,g)(x)=0.
\]
\end{lem}
\begin{proof}
Let $f,g \in \mathcal L^2$. Using the asymptotic behavior of $xw(x)$ for large $x$, see Lemma \ref{lem:asymp w&u}, we find that $f$ and $g$ satisfy
\[
\lim_{k  \rightarrow \infty} s^{-k}  f(zq^{-k})=\lim_{k  \rightarrow \infty} s^{-k}  g(zq^{-k})=0.
\]
From Definition \ref{def:Casorati} and the asymptotic behavior of $u(x)/x$ from Lemma \ref{lem:asymp w&u} we now see that $\lim_{k \rightarrow \infty} D(f,g)(zq^{-k})=0$.
\end{proof}

\subsection{Self-adjointness}
For $f \in F(\R_q)$ we denote
\begin{align*}
f(0^-)&=\lim_{k \rightarrow \infty}f(z_-q^k), & f(0^+)&=\lim_{k \rightarrow \infty} f(z_+q^k), \\
f'(0^-)&=\lim_{k \rightarrow \infty}(D_qf)(z_-q^k), & f'(0^+)&=\lim_{k \rightarrow \infty} (D_qf)(z_+q^k),
\end{align*}
provided that all these limits exist. 
\begin{Def}
We define the subspace $\mathcal D \subset \mathcal L^2$ by
\[
\mathcal D= \left\{ f \in \mathcal L^2 \ | \ Lf \in \mathcal L^2,\ f(0^-)=f(0^+),\ f'(0^-)= f'(0^+) \right\}.
\]
\end{Def}
The domain $\mathcal D$ contains the finitely supported functions in $\mathcal L^2$, hence $\mathcal D$ is dense in $\mathcal L^2$.

\begin{prop}
The operator $(L,\mathcal D)$ is self-adjoint.
\end{prop}
The proposition is proved in the same way as \cite[Prop.2.7]{KS03}. For convenience we repeat the proof here.
\begin{proof}
First we need to show that $(L,\mathcal D)$ is symmetric. Let $f,g \in \mathcal D$. Using the second expression in Definition \ref{def:Casorati} we find
\[
\begin{split}
D(f,\overline{g})(0^-) &= u(0)\Big( (D_qf)(0^-)\overline{g(0^-)} - f(0^-) \overline{ (D_qg)(0^-)} \Big)\\
&= u(0)\Big((D_qf)(0^+)\overline{g(0^+)} - f(0^+) \overline{ (D_qg)(0^+)} \Big) \\
&= D(f,\overline{g})(0^+),
\end{split}
\]
By Proposition \ref{prop:Casorati} and Lemma \ref{lem:lim infty} this leads to
\[
\begin{split}
\langle Lf, g \rangle_{\mathcal L^2}- \langle f,Lg \rangle_{\mathcal L^2} &= \lim_{\substack{l,m \rightarrow \infty \\ k,n \rightarrow -\infty}}\Big( \langle Lf, g \rangle_{k,l;m,n}- \langle f,Lg \rangle_{k,l;m,n} \Big)\\
& =  D(f,\overline{g})(0^+)-D(f,\overline{g})(0^-)=0,
\end{split}
\]
hence $(L,\mathcal D)$ is symmetric with respect to $\langle \cdot, \cdot \rangle_{\mathcal L^2}$.

Now we know that $(L,\mathcal D) \subset (L^*, \mathcal D^*)$, where $(L^*, \mathcal D^*)$ is the adjoint of the operator $(L,\mathcal D)$. Observe that 
\[
L^* = L|_{\mathcal D^*}.
\]
Indeed, let $f$ be a nonzero function with support at only one point $x \in \R_q$ and let $g \in F(\R_q)$, then 
\[
\langle Lf,g \rangle_{\mathcal L^2}=\langle f, Lg \rangle_{\mathcal L^2}.
\]
In particular, for $g \in \mathcal D^*$ we then have $\langle f, Lg \rangle_{\mathcal L^2} = \langle f, L^*g \rangle_{\mathcal L^2}$, so $(Lg)(x)=(L^*g)(x)$. This holds for all $x \in \R_q$, hence $L^*=L|_{\mathcal D^*}$.

Finally we show that $\mathcal D^* \subset \mathcal D$. Let $f \in \mathcal D$ and let $g \in \mathcal D^*$. Using Proposition \ref{prop:Casorati} and Lemma \ref{lem:lim infty}, 
\[
D(f,\overline{g})(0^-)-D(f,\overline{g})(0^+) = \langle Lf,g \rangle_{\mathcal L^2} - \langle f,L^*g \rangle_{\mathcal L^2}=0.
\]
Since this holds for all $f \in \mathcal D$, we find that the limits $g(0^-)$, $g(0^+)$, $g'(0^-)$ and $g'(0^+)$ exist, and
\[
\begin{split}
g(0^-) =  g(0^+), \qquad g'(0^-)&=  g'(0^+),
\end{split}
\]
hence $g \in \mathcal D$, which proves the proposition.
\end{proof}
\begin{rem}
Let $f \in \mathcal D$ and let $\al$ be a complex number with $|\al|=1$. We define
\[
\tilde f(x) = 
\begin{cases}
f(x),& x \in \R_q^-,\\
\overline{\al} f(x), & x \in \R_q^+,
\end{cases}
\]
then it is easy to verify that $\tilde f(0^-)=\al \tilde f(0^+)$ and $\tilde f'(0^-) = \al \tilde f'(0^+)$. So we have a family of dense domains
\[
\mathcal D_\al = \left\{ f \in \mathcal L^2 \ | \ Lf \in \mathcal L^2,\ f(0^-)=\al f(0^+),\ f'(0^-)= \al f'(0^+) \right\}, 
\]
such that $(L,\mathcal D_\al)$ is self-adjoint. Without loss of generality we may work with the dense domain $\mathcal D=\mathcal D_1$.
\end{rem}

\section{Eigenfunctions} \label{sec3}
In this section we study eigenfunctions of the second order difference operator $L$.

\subsection{Spaces of eigenfunctions}
For $\mu \in \C$ we introduce the spaces
\[
\begin{split}
V_\mu^-&=\{ f: \R_q^- \rightarrow \C \ | \ Lf=\mu f \},\\
V_\mu^+&=\{ f: \R_q^+ \rightarrow \C \ | \ Lf =\mu f \},\\
V_\mu &= \{ f: \R_q \rightarrow \C \ | \ Lf=\mu f, \ f(0^-)=f(0^+), \ f'(0^-)=f'(0^+) \}
\end{split}
\]
\begin{lem} \label{lem:dim V}
Let $\mu \in \C$.

\textnormal{(a)} For $f,g \in V_\mu^\pm $ the Casorati determinant $D(f,g)$ is constant on $\R_q^\pm$. 

\textnormal{(b)} For $f,g \in V_\mu$ the Casorati determinant $D(f,g)$ is constant on $\R_q$.

\textnormal{(c)} $\dim V_\mu^\pm = 2$. 

\textnormal{(d)} $\dim V_\mu \leq 2$.
\end{lem}
\begin{proof}
For (a) let $f,g \in F(\R_q)$. From the proof of Proposition \ref{prop:Casorati} we have
\[
\Big((Lf)(x)g(x) - f(x) (Lg)(x)\Big) (1-q)xw(x) = D(f,g)(x/q)- D(f,g)(x).
\]
Now if $f$ and $g$ satisfy $(Lf)(x)=\mu f(x)$ and $(Lg)(x)=\mu g(x)$, we find $D(f,g)(x/q)= D(f,g)(x)$, hence $D(f,g)$ is constant on $\R_q^+$ and $\R_q^-$.

Let $f,g \in V_\mu$. Statement (b) follows from (a) and the fact that $D(f,g)(0^-) =  D(f,g)(0^+)$. 

For (c) we write $f(zq^k)= f_k$, then we see that $Lf =\mu f$ gives a recurrence relation of the form $\al_k f_{k+1} +\be_k f_k + \ga_k f_{k-1} = \mu f_k$, with $\al_k, \ga_k \neq 0$ for all $k \in \Z$. Solutions of such a recurrence relation are uniquely determined by specifying $f_k$ at two different points $k=l$ and $k=m$. So there are two independent solutions, which means that $\dim V_\mu^\pm=2$.

Finally, suppose that $f_1,f_2 \in V_\mu$ are such that the restrictions $f_i^{res}=f_i|_{\R_q^+}$ are linearly independent in $V_\mu^+$. By (a) the Casorati determinant $D(f_1^{res},f_2^{res})(x)$ is nonzero and constant on $\R_q^+$, hence $D(f_1,f_2)$ is nonzero and constant on $\R_q$. Therefore $f_1$ and $f_2$ are linearly independent. Now choose a function $f_3 \in V_\mu$. Since $\dim V_\mu^+=2$, we have $f_3^{res} = \al f_1^{res} + \be f_2^{res}$ for some constants $\al,\be \in \C$. This shows that 
\[
\begin{split}
D(f_3,f_1) = D(f_3^{res},f_1^{res}) &= \be D(f_2^{res},f_1^{res}) = \be D(f_2,f_1), \\
D(f_3,f_2) = D(f_3^{res},f_2^{res}) &= \al D(f_1^{res},f_2^{res}) = \al D(f_1,f_2),
\end{split}
\]
hence $f_3 = \al f_1 +\be f_2$. So $\dim V_\mu\leq2$.
\end{proof}

\subsection{Big $q$-Jacobi functions} 
Let $P_{\mathrm{gen}}$ be the dense subset of $P$ given by
\[
P_{\mathrm{gen}} = \{ (a,b,c,d) \in P \ |\ c \neq d, \ c/a,c/b,d/a,d/b,cd/ab \not\in q^\Z\}.
\]
From here on we assume that $(a,b,c,d) \in P_{\mathrm{gen}}$, unless stated otherwise. 

The difference operator $L$ is equivalent to the difference operator studied in \cite{KS03}. To see this set
\begin{equation} \label{eq:param KS}
X=-\frac{ax}{q},\quad A=s, \quad B=\frac{qc}{sa},\quad C= \frac{qd}{sa},
\end{equation}
where the capitals stand for the parameters in \cite{KS03}. We have the following eigenfunction, which is called a big $q$-Jacobi function,
\begin{equation} \label{def:phi}
\varphi_\ga(x) = \varphi_\ga(x;a,b,c,d|q) =  \ph{3}{2}{q/ax, s\ga, s/\ga}{cq/a, dq/a}{q,bx}, \qquad |x| <  \frac{1}{|b|} .
\end{equation}
If $|x|<|q/b|$, the function $\varphi_\ga$ is a solution of the eigenvalue equation
\begin{equation} \label{eq:eigenvalue equation}
(Lf)(x) = \mu(\ga)f(x), \qquad \mu(\ga)=\ga+\ga^{-1},
\end{equation}
where $\ga \in \C^*$ and $x \in \R_q$. This can be obtained from \cite{KS03}, or directly from the contiguous relation \cite[(2.10)]{GIM96}. For a function $f$ depending on the parameters $a,b,c,d$, $f = f(\cdot;a,b,c,d)$, we write
\[
f^\dagger=f^\dagger(\cdot;a,b,c,d)= f(\cdot;b,a,c,d).
\]
Clearly, we have $(f^\dagger)^\dagger =f$. Since $L_{a,b,c,d}=L_{b,a,c,d}$, cf.~ Remark \ref{rem:symmetry}(a), it is immediately clear that $\varphi_\ga^\dagger$ is also a solution for the eigenvalue equation \eqref{eq:eigenvalue equation}. If $a=\overline{b}$, we have $\varphi_\ga^\dagger(x) = \overline{\varphi_{\overline{\ga}}(x)}$. The symmetry $L_{a,b,c,d}=L_{a,b,d,c}$ does not give rise to different eigenfunctions.

So far the functions $\varphi_\ga(x)$ and $\varphi_\gamma^\dagger(x)$ are defined for small $x \in \R_q$. Using the eigenvalue equation \eqref{eq:eigenvalue equation} the functions $\varphi_\ga$ and $\varphi_\ga^\dagger$ can uniquely be extended to functions on whole $\R_q$ (that we also denote by $\varphi_\ga$ and $\varphi_\ga^\dagger$) that also satisfy \eqref{eq:eigenvalue equation}. Later on we give explicit expressions for the functions $\varphi_\ga(x)$ and $\varphi_\ga^\dagger(x)$ for $|x|>q/|b|$.

First we establish the $q$-differentiability at the origin of the functions $\varphi_\ga$ and $\varphi_\ga^\dagger$.
\begin{prop} \label{prop:diff origin}
The functions $\varphi_\ga$ and $\varphi_\ga^\dagger$ are continuous differentiable at the origin. At $x=0$ we have 
\[
\begin{split}
\varphi_\ga(0;a,b,c,d|q) &= \ph{2}{2}{\ga\sqrt{cdq/ab}, \sqrt{cdq/ab}/\ga}{cq/a, dq/a}{q,\frac{bq}{a} }, \\ \varphi_\ga'(0;a,b,c,d|q)&=\frac{b(1-s\ga)(1-s/\ga)}{(1-q)(1-cq/a)(1-dq/a)} \varphi_\ga(0;aq^{-\hf},bq^{-\hf},cq^\hf, dq^\hf|q).
\end{split}
\]
\end{prop}
\begin{proof}
The expression for $\varphi_\ga(0)$ follows from letting $x \rightarrow 0$ in \eqref{def:phi}.
If $|x|$ is small enough, we find from the explicit expression \eqref{def:phi} for $\varphi_\ga$,
\[
\begin{split}
\varphi_\ga(x)-\varphi_\ga(qx) = &\sum_{n=1}^\infty \frac{ (s\ga,s/\ga;q)_n }{ (q,cq/a,dq/a;q)_n } (bx)^n \Big[ (q/ax;q)_n - (1/ax;q)_n q^n \Big] \\
=& \sum_{n=1}^\infty \frac{ (s\ga,s/\ga;q)_n }{ (q,cq/a,dq/a;q)_n } (bx)^n (q/ax;q)_{n-1} (1-q^n)\\
=& \frac{ bx (1-s\ga)(1-s/\ga) }{ (1-cq/a)(1-dq/a)} \ph{3}{2}{ q/ax, sq \ga, sq/\ga}{cq^2/a, dq^2/a}{q,bx}.
\end{split}
\]
Now it follows that the $q$-derivative of $\varphi_\ga$ is given by
\[
(D_q \varphi_\ga)(x) = \frac{b(1-s\ga)(1-s/\ga)}{(1-q)(1-cq/a)(1-dq/a)} \varphi_\ga(xq^\hf;a q^{-\hf},bq^{-\hf}, cq^{\hf}, dq^\hf  |q).
\]
Letting $x \rightarrow 0$ gives the result.
\end{proof}

\subsection{Asymptotic solutions}
We define the set of regular spectral values
\[
\mathcal S_{\mathrm{reg}} = \C^*\setminus \{\pm q^{\hf k} \ |\ k \in \Z \}.
\]
For $\ga \in \mathcal S_{\mathrm{reg}}\cup \{ \pm 1\}$ another solution for the eigenvalue equation \eqref{eq:eigenvalue equation} is the function
\begin{equation} \label{def:Phi}
\begin{split}
\Phi_\ga(x)&=\Phi_\ga(x;a,b,c,d|q) \\
&= (s\ga)^k \frac{ (q/bx, q^2\ga/asx;q)_\infty }{ (q/cx, q/dx ; q)_\infty } \ph{3}{2}{ q \ga/s, cq\ga/sa, dq\ga/sa}{ q^2 \ga /asx, q\ga^2}{q, \frac{ q}{bx} }, \quad |x| > \frac{q}{|b|},
\end{split}
\end{equation}
where $x=zq^{-k}$, see \cite{KS03}. For $x \rightarrow \pm \infty$ we have
\begin{equation} \label{eq:asymp Phi}
\Phi_\ga(zq^{-k})= (s\ga)^k(1+\mathcal O(q^k)), \quad k \rightarrow \infty.
\end{equation}
Clearly $\Phi_{\ga^{-1}}$, $\Phi_{\ga}^\dagger$ and $\Phi_{\ga^{-1}}^\dagger$ are also solutions to \eqref{eq:eigenvalue equation}. We remark that it follows from applying the transformation \cite[(III.9)]{GR04} for $_3\phi_2$-series, that
\[
\Phi_\ga(x) = (s\ga)^k \frac{ (q/ax, q^2\ga/bsx;q)_\infty }{ (q/cx, q/dx;q)_\infty } \ph{3}{2}{ q \ga/s, cq\ga/sb, dq\ga/sb}{ q^2 \ga /bsx, q\ga^2}{q, \frac{q}{ax} }, \quad x=zq^{-k}.
\]
So we see that $\Phi_\ga=\Phi_{\ga}^\dagger$, and if $\ga \in \R$ we see that $\Phi_\ga$ is real-valued. Using the eigenvalue equation $L\Phi_{\ga}=\mu(\ga)\Phi_\ga$ we can extend $\Phi_\ga$ to single-valued functions $\Phi^+_\ga$ on $\R_q^+$, respectively $\Phi^-_\ga$ on $\R_q^-$. We call $\Phi_\ga^+$ and $\Phi_\ga^-$ the asymptotic solutions of $Lf=\mu(\ga)f$ on $\R_q^+$, respectively $\R_q^-$. The following lemma shows that $\Phi_\ga^\pm$ and $\Phi_{\ga^{-1}}^\pm$ are linear independent, hence they form linear bases for the eigenspaces $V_{\mu(\ga)}^\pm$.
\begin{lem} \label{lem:det Phi}
For $\ga \in \mathcal S_{\mathrm{reg}}$ we have 
\[
D(\Phi_\ga^\pm,\Phi_{1/\ga}^\pm)(z_{\pm}q^{-k}) =  (\ga-1/\ga)K_{z_\pm}.
\]
\end{lem}
\begin{proof}
Since $\Phi_\ga^+$ lies in $V^+_{\mu(\ga)}$, the Casorati determinant $D(\Phi_\ga^+,\Phi_{1/\ga}^+)$ is constant on $\R_q^+$, so we can find the Casorati determinant by taking the limit $x \rightarrow  \infty$. From Lemma \ref{lem:asymp w&u} we find
\[
\lim_{k \rightarrow \infty} \frac{s^{2k-1} u(z_+q^{-k})}{z_+q^{-k}(1-q)} = K_{z_+},
\]
and then it follows from the first expression in Definition \ref{def:Casorati} and \eqref{eq:asymp Phi} that
\[
\lim_{k \rightarrow \infty} D(\Phi_\ga^+,\Phi_{1/\ga}^+)(z_+q^{-k})=  (\ga-1/\ga)K_{z_+}. 
\]
The proof for $\Phi^-_\ga$ is similar.
\end{proof}
Now we can expand the functions $\varphi_\ga$ and $\varphi_\ga^\dagger$ on $V_{\mu(\ga)}^\pm$ in terms of $\Phi_{\ga^{\pm 1}}^\pm$. The expansion of $\varphi_\ga$ in terms of $\Phi_{\ga^{\pm1}}^+$ (respectively $\Phi_{\ga^{\pm1}}^-$) gives an explicit expression for $\varphi_\ga$ for $x>q/|b|$ (respectively $x<-q/|b|$). For $\ga \in \mathcal S_{\mathrm{reg}}$ we define a function $c_z(\ga)$ by
\[
c_z(\ga) = c_z(\ga;a,b,c,d|q) = \frac{ (s/\ga, cq/as\ga, dq/as\ga;q)_\infty \te( bsz\ga ) } { (cq/a, dq/a, 1/\ga^2;q)_\infty \te(bz) }.
\]
The desired expansion uses the $c$-function, see \cite[Prop.4.4]{KS03} with parameters as in \eqref{eq:param KS} and $Z = za/q$, or use the three-term transformation for $_3\phi_2$-functions \cite[(III.33)]{GR04}.
\begin{prop} \label{prop:c function exp}
For $\ga \in \mathcal S_{\mathrm{reg}}$ and $x = z_{\pm}q^k \in \R_q$,
\[
\begin{split}
\varphi_\ga(x) &= c_{z_\pm}(\ga) \Phi_\ga^\pm(x) + c_{z_\pm}(\ga^{-1}) \Phi_{\ga^{-1}}^\pm(x),\\
\varphi_\ga^\dagger(x) &= c^\dagger_{z_\pm}(\ga) \Phi_\ga^\pm(x) + c^\dagger_{z_\pm}(\ga^{-1}) \Phi_{\ga^{-1}}^\pm(x).
\end{split}
\]
\end{prop}

The spaces $V^\pm_2$ and $V^\pm_{-2}$ are 2-dimensional by Lemma \ref{lem:dim V}, but they are clearly not spanned by $\Phi_\ga^\pm$ and $\Phi_{1/\ga}^\pm$, since $\ga=\pm 1$ here. In the following lemma we give linear bases for the spaces $V^\pm_2$ and $V^\pm_{-2}$ that will be useful later on.
\begin{lem} \label{lem:ga = 1}
For $\ga=1$ or $\ga=-1$, the functions $\Phi_\ga^\pm$ and $\left.\dfrac{d\Phi_{\ga'}^\pm}{d\ga'}\right|_{\ga'=\ga}$ form a linear basis for the spaces $V^\pm_2$, respectively $V^\pm_{-2}$. 
\end{lem}
\begin{proof}
Differentiating the equation $L\Phi_\ga^+=\mu(\ga)\Phi_\ga^+$ to $\ga$, and setting $\ga=\pm 1$ shows that $\frac{d\Phi_\ga^+}{d\ga}|_{\ga=\pm 1}$ is an eigenfunction of $L$ for eigenvalue $\pm 2$. From the asymptotic behavior \eqref{eq:asymp Phi} of $\Phi_\ga^+$ we find
\[
\frac{d\Phi_\ga^+}{d\ga}(zq^{-k}) = sk (s\ga)^{k-1}\Big(1+\mathcal O(q^k)\Big), \quad k \rightarrow \infty,
\]
and then using Lemma \ref{lem:asymp w&u} it follows that
\[
D\Big(\Phi_\ga^+, \frac{d\Phi_\ga^+}{d\ga}\Big)(x)= \ga^{2k-2} K_{z_+}, \qquad x \in \R_q^+.
\]
For $\ga =\pm1 $ we see that $D(\Phi_{\pm 1}^+,\left.\frac{d\Phi_\ga^+}{d\ga}\right|_{\ga=\pm 1})= K_{z_+} \neq 0$. This proves the lemma for $\Phi_{\ga}^+$. For $\Phi_\ga^-$ the proof is the same.
\end{proof}

\subsection{A basis for $V_\mu$}
We are going to show that, under certain conditions on $\ga$, the solutions $\varphi_\ga$ and $\varphi_\ga^\dagger$ form a linear basis for $V_{\mu(\ga)}$. We do this by computing the Casorati determinant $D(\varphi_\ga,\varphi_{\ga}^\dagger)$. 
\begin{lem} \label{lem:det phi}
For $x \in \R_q$ and $\ga\in \C^*$ we have
\[
D(\varphi_\ga, \varphi_\ga^\dagger)(x) = \frac{(1-q)q}{as} \frac{ (s\ga,s/\ga;q)_\infty\te(a/b)} {(cq/a,cq/b,dq/a,dq/b;q)_\infty}.
\]
\end{lem}
\begin{proof}
Let $\ga \in \mathcal S_{\mathrm{reg}}$. From Proposition \ref{prop:diff origin} we know that $\varphi_\ga, \varphi_\ga^\dagger \in V_{\mu(\ga)}$, hence by Lemma \ref{lem:dim V} the Casorati determinant $D(\varphi_\ga, \varphi_\ga^\dagger)$ is constant on $\R_q$. In order to calculate $D(\varphi_\ga, \varphi_\ga^\dagger)(x)$ we use the $c$-function expansions from Proposition \ref{prop:c function exp}, 
\[
D(\varphi_\ga, \varphi_\ga^\dagger)(x) = \sum_{\epsilon,\eta\in\{-1,1\}}c_{z_+}(\ga^{\epsilon})c_{z_+}^\dagger(\ga^\eta) D(\Phi_{\ga^\epsilon}^+,\Phi_{\ga^\eta}^+)(x), \qquad x = z_+q^{-k}.
\]
We apply Lemma \ref{lem:det Phi}, then
\[
D(\varphi_\ga, \varphi_\ga^\dagger)(z_+q^{-k}) = (\ga-1/\ga) K_{z_+} \Big( c_{z_+}(\ga)c^\dagger_{z_+}(1/\ga) - c_{z_+}(1/\ga)c^\dagger_{z_+}(\ga) \Big).
\]
Using $cq/as= bs/d$ and $dq/as=bs/c$, we find
\[
\begin{split}
 c_{z_+}(\ga)&c^\dagger_{z_+}(1/\ga) - c_{z_+}(1/\ga)c^\dagger_{z_+}(\ga)  = \\
&\frac{ (s\ga,s/\ga;q)_\infty }{ (\ga^2,1/\ga^2;q)_\infty  (cq/a,cq/b,dq/a,dq/b;q)_\infty \te(az_+,bz_+) } \\ 
\times&  \Big(\te(q/bsz_+\ga, q\ga/asz_+, cq\ga/bs, cq/as\ga) -  \te(q\ga/bsz_+, q/asz_+\ga, cq/bs\ga, cq\ga/as) \Big)
\end{split}
\]
Now we use the $\te$-product identity \eqref{eq:theta prod} with
\begin{align*}
x&= \frac{qe^{\hf i\kappa}}{bs} \sqrt{\frac{|c|}{z_+}},& v&=\ga e^{\hf i\kappa} \sqrt{|c|z_+}, \\ y&=\frac{qe^{\hf i\kappa}}{as} \sqrt{\frac{|c|}{z_+}},& w&=\ga e^{-\hf i\kappa}\sqrt{\frac{1}{|c|z_+}},
\end{align*}
where $c=|c|e^{i\kappa}$, then we obtain
\[
c_{z_+}(\ga)c^\dagger_{z_+}(1/\ga) - c_{z_+}(1/\ga)c^\dagger_{z_+}(\ga)
= \frac{ q}{asz_+(\ga-1/\ga)}  \frac{(s\ga,s/\ga;q)_\infty \te(a/b,cz_+,dz_+) }{(cq/a,cq/b,dq/a,dq/b;q)_\infty \te(az_+,bz_+)} .
\]
With the explicit expression for $K_{z_+}$ we find the desired result for $\ga \in \mathcal S_{\mathrm{reg}}$. By continuity in $\ga$ the result holds for all $\ga \in \C^*$.
\end{proof}
Let $\mathcal S_{\mathrm{pol}}$ be the set of zeros of $\ga \mapsto D(\varphi_\ga,\varphi_\ga^\dagger)$, i.e.,
\[
\mathcal S_{\mathrm{pol}} = \left\{ sq^k \ | \ k \in \Z_{\geq 0} \right\} \cup \left\{ s^{-1}q^{-k} \ | \ k \in \Z_{\geq 0} \right\}.
\]
\begin{prop} \label{prop:basis Vmu}
Let $\ga \in \C^*\setminus \mathcal S_{\mathrm{pol}}$, then $\dim V_{\mu(\ga)}=2$ and the set $\{\varphi_\ga,\varphi_\ga^\dagger\}$ is a linear basis of $V_{\mu(\ga)}$. 
\end{prop}
\begin{proof}
From Lemma \ref{lem:det phi} it follows that $\varphi_\ga$ and $\varphi_\ga^\dagger$ are linearly independent if $\ga \not\in \mathcal S_{\mathrm{pol}}$. Since both functions are continuously differentiable at the origin, see Proposition \ref{prop:diff origin}, and since $\dim V_{\mu(\ga)}\leq 2$ by Lemma \ref{lem:dim V}, we have $\dim V_{\mu(\ga)}=2$, and $\varphi_\ga$ and $\varphi_\ga^\dagger$ form a linear basis for $V_{\mu(\ga)}$. 
\end{proof}
\begin{cor} \label{cor:extensions}
For $\ga \in \C^*\setminus \mathcal S_{\mathrm{pol}}$ every function in $V_{\mu(\ga)}^+$ (respectively $V_{\mu(\ga)}^-$) has a unique extension to $V_{\mu(\ga)}$.
\end{cor}
\begin{proof}
Fix a $\ga \in \C^*\setminus \mathcal S_{\mathrm{pol}}$ and denote $\mu = \mu(\ga)$. We consider the restriction map $res:V_\mu \rightarrow V_\mu^+$ defined by $f^{res} = f|_{\R_q^+}$. Let $f$ and $g$ be linearly independent in $V_\mu$. As in the proof of Lemma \ref{lem:dim V} it follows that $f^{res}$ and $g^{res}$ are linearly independent in $V_\mu^+$. Since $\dim V_\mu = \dim V_\mu^+ =2$ the map ${res}$ is a linear isomorphism. In a similar way a linear isomorphism between $V_\mu$ and $V_\mu^-$ can be constructed.
\end{proof}

For $\ga \in \mathcal S_{\mathrm{pol}}$ the big $q$-Jacobi functions $\varphi_\ga$ and $\varphi_\ga^\dagger$ are actually multiples of big $q$-Jacobi polynomials, see \cite[Prop.5.3]{KS03}. The big $q$-Jacobi polynomials, see \cite{AA85}, \cite{KS98}, are defined by
\[
P_k(x;\al,\be,\de;q) = \ph{3}{2}{q^{-k}, \al \be q^{k+1}, x}{\al q, \de q}{q,q}, \qquad k \in \Z_{\geq 0}.
\]

\begin{lem} \label{lem:phi pol}
Let $\ga_k=sq^k \in \mathcal S_{\mathrm{pol}}$ or $\ga_k=s^{-1}q^{-k} \in \mathcal S_{\mathrm{pol}}$, then
\[
\varphi_{\ga_k}(x) = \frac{ (cq/b, dq/b;q)_k }{ (cq/a, dq/a;q)_k } \left( \frac{b}{a} \right)^k \varphi_{\ga_k}^\dagger(x), 
\]
and 
\[
\varphi_{\ga_k}^\dagger(x) =  q^{-\hf k (k+1)} \left( -\frac{ a}{c} \right)^k \frac{ (cq/a;q)_k }{ (dq/b;q)_k } P_k(cx;c/b, d/a, c/a;q),
\]
for $x \in \R_q$. 
\end{lem}
\begin{proof}
Let $k \in \Z_{\geq 0}$ and $\ga_k= sq^k$. From Lemma \ref{lem:det phi} we see that the Casorati determinant $D(\varphi_{\ga_k},\varphi_{\ga_k}^*)(x)$, $x \in \R_q$, is equal to zero, hence $\varphi_{\ga_k}(x) = C_k \varphi_{\ga_k}^\dagger(x)$, for some constant $C_k$ independent of $x$. To find the constant $C_k$ we use Proposition \ref{prop:c function exp}. We have $c_{z}(\ga_k)=0$ and $c_{z}^\dagger(\ga_k)=0$, hence 
\[
\varphi_{\ga_k}(x) = c_{z_+}(1/\ga_k)\Phi_{1/\ga_k}^+(x) = \frac{c_{z_+}(1/\ga_k)}{c_{z_+}^\dagger(1/\ga_k)}\varphi_{\ga_k}^\dagger(x), \qquad x=z_+q^n \in \R_q.
\]
Using $\te(q^k x) = (-x)^{-k} q^{-\hf k (k-1)} \te(x)$ we find
\[
C_k=\frac{c_{z_+}(1/\ga_k)}{c_{z_+}^\dagger(1/\ga_k)} = \frac{ (cq/b, dq/b;q)_k }{ (cq/a, dq/a;q)_k } \left( \frac{b}{a} \right)^k.
\]
Since $\varphi_\ga = \varphi_{1/\ga}$ the result also holds for $\ga_k=s^{-1}q^{-k}$, $k \in \Z_{\geq 0}$.

Finally, writing out $\varphi_{\ga_k}^\dagger(x)$ explicitly as a $_3\phi_2$-series using \eqref{def:phi} and then applying the transformation formula \cite[(III.13)]{GR04} shows that $\varphi_{\ga_k}^\dagger(x)$ is a multiple of a big $q$-Jacobi polynomial in the variable $cx$.
\end{proof}

\subsection{Extensions of the asymptotic solutions}
By Corollary \ref{cor:extensions} the asymptotic solutions $\Phi_\ga^+ \in V_{\mu(\ga)}^+$ and $\Phi_\ga^- \in V_{\mu(\ga)}^-$ have unique extensions to $V_{\mu(\ga)}$, provided that $\ga \in \C^*\setminus \mathcal S_{\mathrm{pol}}$. We denote these extensions again by $\Phi_\ga^+$ and $\Phi_\ga^-$. Propositions \ref{prop:basis Vmu} and \ref{prop:c function exp} enable us to expand $\Phi_{\ga^{\pm 1}}^\pm$ in terms of the basis $\{ \varphi_\ga, \varphi_\ga^\dagger\}$ of $V_{\mu(\ga)}$. 
\begin{prop} \label{prop:Phi+-}
For $x \in \R_q$ and $\ga \in \mathcal S_{\mathrm{reg}}\setminus \{sq^k \ | \ k \in \Z_{\geq 0} \}$
\[
\begin{split}
\Phi^+_\ga(x) &= d_{z_+}(\ga)\varphi_\ga(x)+ d^\dagger_{z_+}(\ga)\varphi_\ga^\dagger(x),\\
\Phi^-_\ga(x) &= d_{z_-}(\ga)\varphi_\ga(x)+ d^\dagger_{z_-}(\ga)\varphi_\ga^\dagger(x),
\end{split}
\]
where
\[
d_z(\ga)=d_z(\ga;a,b,c,d|q) = \frac{  (cq/a, dq/a;q)_\infty \te(bz) }{  \te(a/b,cz,dz) } \ \frac{ (cq\ga/sb, dq\ga/sb;q)_\infty \te(asz/q\ga) }{ (q\ga^2,s/\ga;q)_\infty}.
\]
\end{prop}
\begin{proof}
Let $\ga \in \mathcal S_{\mathrm{reg}}\setminus \mathcal S_{\mathrm{pol}}$. By Proposition \ref{prop:basis Vmu} we may expand 
\[
\Phi_\ga^\pm(x) = d_{z_\pm}(\ga) \varphi_\ga(x) + e_{z_\pm}(\ga) \varphi_\ga^\dagger(x), \qquad x \in \R_q,
\]
for some coefficients $d_{z}(\ga)$ and $e_z(\ga)$ independent of $x$. In order to compute the coefficients $d_z(\ga)$ and $e_z(\ga)$ we observe that it follows from $\Phi_\ga^{\pm \dagger} = \Phi_\ga^\pm$ that $e_z(\ga) = d_z^\dagger(\ga)$. To compute $d_z(\ga)$ we use
\[
d_{z_{\pm}}(\ga) = \frac{ D(\Phi_\ga^\pm , \varphi_\ga^\dagger)(x) }{ D(\varphi_\ga, \varphi_\ga^\dagger)(x) }.
\]
From the $c$-function expansion, Proposition \ref{prop:c function exp}, we find
\[
D(\Phi_\ga^\pm, \varphi_\ga^\dagger)(x) = c^\dagger_{z_\pm}(1/\ga)D(\Phi_\ga^\pm,\Phi_{1/\ga}^\pm)(x), \qquad x \in \R_q,
\]
and then it follows from Lemmas \ref{lem:det Phi} and \ref{lem:det phi} that
\[
\begin{split}
d_z(\ga) &= \frac{(\ga-1/\ga)K_z \,c^\dagger_z(1/\ga) }{D(\varphi_\ga,\varphi_\ga^\dagger)}  \\
&= \frac{ asz (cq/a, dq/a;q)_\infty \te(bz) }{ q\, \te(a/b,cz,dz) } \ \frac{(\ga-1/\ga) (cq\ga/bs, dq\ga/bs;q)_\infty \te(asz/\ga) }{ (\ga^2,s/\ga;q)_\infty }.
\end{split}
\]
Here we also used the explicit expression for $K_z$ from Lemma \ref{lem:asymp w&u}. This is the desired result for $\ga \in \mathcal S_{\mathrm{reg}}\setminus \mathcal S_{\mathrm{pol}}$. By continuity in $\ga$ it holds also for $\ga \in \{ s^{-1}q^{-k} \ | \ k \in \Z\}$.
\end{proof}
For $\ga=s^{-1}q^{-k}$, $k \in \Z_{\geq 0}$, the Casorati determinant $D(\varphi_\ga,\varphi_\ga^\dagger)$ is equal to zero, hence $\varphi_\ga$ is a multiple of $\varphi_\ga^\dagger$. In this case Proposition \ref{prop:Phi+-} states that $\Phi_\ga^\pm$ is also a multiple of $\varphi_\ga^\dagger$.
\begin{cor} \label{cor:Phi pol}
For $\ga_k = s^{-1}q^{-k}$, $k \in \Z_{\geq 0}$,
\[
\Phi_{\ga_k}^\pm(x) = q^{\hf k(k-1)}\left(\frac{-1}{az_{\pm}} \right)^{k} \frac{ (cq/b,dq/b;q)_k }{ (s^2;q)_k  } \varphi_{\ga_k}^\dagger(x).
\]
\end{cor}
\begin{proof}
Using Proposition \ref{prop:Phi+-} and Lemma \ref{lem:phi pol}, we find
\[
\Phi_{\ga_k}^+(x) = \frac{b}{aq}\left(\frac{b}{qcdz_+} \right)^{k-1} \frac{ (cq/b,dq/b;q)_k\ F_{z_+}}{ (q^{1-k}/s^2;q)_k \te(b/a,cz_+, dz_+, s^2) } \varphi_{\ga_k}^\dagger(x),
\]
with 
\[
F_{z_+} = \te(cq/a, dq/a, q/bz_+, cdqz_+/b) - \te( cq/b, dq/b, q/az_+, cdqz_+/a).
\]
Applying the $\te$-product identity \eqref{eq:theta prod} with
\[
x=\frac{qe^{i(\kappa+\de)}\sqrt{|cd|}}{b}, \quad y=\frac{qe^{i(\kappa+\de)}\sqrt{|cd|}}{a}, \quad v= \frac{e^{-i(\kappa+\de)}}{z_+\sqrt{|cd|}}, \quad w= e^{i(\kappa-\de)}\sqrt{\left|\frac{c}{d}\right|},
\]
where $c=|c|e^{i\kappa}$ and $d=|d|e^{i\de}$, we obtain
\[
F_{z_+} = \frac{ cdqz_+}{a} \te(cdq^2/ab, a/b, 1/dz_+, 1/cz_+).
\]
Applying $(q^{1-k}/y;q)_k=(-y)^{-k}q^{-\hf k(k-1)} (y;q)_k$, identities for $\te$-functions, and  $s^2=cdq/ab$, the result follows for $\Phi^+_{\ga_k}$. Replacing $z_+$ by $z_-$ gives the result for $\Phi^-_{\ga_k}$.
\end{proof}

In the expansion of $\Phi_\ga^\pm$ in Proposition \ref{prop:Phi+-} we have assumed that $\ga \not\in \{sq^k \ | \ k \in \Z\}$. At first sight it seems that the functions $\Phi_\ga^\pm(x)$, considered as functions of $\ga$ and with $x \in \R_q$ fixed, have simple poles at the points $\ga= sq^k$, $k \in \Z_{\geq 0}$, which are the poles of the function $d_z(\ga)$. It turns out that the functions $\Phi_\ga^\pm(x)$ are actually analytic at these points. 
\begin{prop} \label{prop:Phi pol}
For a given $x \in \R_q$ the functions $\ga \mapsto \Phi_\ga^\pm(x)$ are analytic on $\mathcal S_{\mathrm{reg}}$. In particular, for $\ga_k=sq^k$, $k \in \Z_{\geq 0}$,
\[
\begin{split}
\Phi_{\ga_k}^\pm(x) &= \left.\Res{\ga=\ga_k}\big( d_{z_\pm}(\ga)\big)\,\frac{ d}{d\ga}\varphi_\ga(x)\right|_{\ga=\ga_k} + \left. \Res{\ga=\ga_k}\big( d_{z_{\pm}}^\dagger(\ga)\big)\,\frac{ d}{d\ga}\varphi_\ga^\dagger(x)\right|_{\ga=\ga_k} + \tilde d_{z_\pm}(\ga_k)\varphi_{\ga_k}^\dagger(x),
\end{split}
\]
where 
\[
\tilde d_{z}(\ga_k) = \lim_{\ga \rightarrow \ga_k} \Bigg(\frac{ (cq/b, dq/b;q)_k }{ (cq/a, dq/a;q)_k } \left( \frac{b}{a} \right)^k d_{z}(\ga) + d_{z}^\dagger(\ga) \Bigg). 
\]
\end{prop}
\begin{proof}
The expansion from Proposition \ref{prop:Phi+-} shows that the functions $\ga \mapsto \Phi_\ga^\pm(x)$, for a given $x \in \R_q$, are analytic functions on $\mathcal S_{\mathrm{reg}} \setminus \{ sq^k \ | \ k \in \Z\}$. So we only have to consider the functions $\Phi_\ga^\pm(x)$ at the points $\ga_k = sq^k$, $k \in \Z_{\geq 0}$.

Fix a $k \in \Z_{\geq 0}$ and a $x \in \R_q$. The function $\ga \mapsto d_z(\ga)$ has a simple pole at $\ga=\ga_k$ coming from the zero of the infinite product $(s/\ga;q)_\infty$, and the functions $\ga \mapsto \varphi_\ga(x)$ and $\ga \mapsto \varphi_\ga^\dagger(x)$ are analytic at $\ga=\ga_k$. From Proposition \ref{prop:Phi+-} and Lemma \ref{lem:phi pol} it follows that
\[
\begin{split}
\Phi_\ga^+(x) =& (\ga-\ga_k) d_{z_+}(\ga) \frac{ \varphi_\ga(x) - \varphi_{\ga_k}(x)}{\ga-\ga_k} + (\ga-\ga_k) d_{z_+}^\dagger(\ga) \frac{\varphi_\ga^\dagger(x) - \varphi_{\ga_k}^\dagger(x)}{\ga-\ga_k}\\
& + \Bigg(\frac{ (cq/b, dq/b;q)_k }{ (cq/a, dq/a;q)_k } \left( \frac{b}{a} \right)^k d_{z_+}(\ga) + d_{z_+}^\dagger(\ga) \Bigg) \varphi_{\ga_k}^\dagger(x).
\end{split}
\]
We see that the limit $\lim_{\ga \rightarrow \ga_k} \Phi_\ga^+(x)$ exists if $\tilde d_{z_+}(\ga_k)$, as defined in the proposition, exists. Let us define
\[
\hat d_z(\ga) = (s/\ga;q)_\infty d_z(\ga),
\]
then $\ga \mapsto \hat d_z(\ga)$ is regular at $\ga=\ga_k$. By a straightforward computation we obtain
\[
\frac{ (cq/b, dq/b;q)_k }{ (cq/a, dq/a;q)_k } \left( \frac{b}{a} \right)^k \hat d_z(\ga_k) + \hat d_z^\dagger(\ga_k) =0,
\]
and then it follows that $\tilde d_z(\ga_k)$ exists.
\end{proof}

We now have the following properties of the functions $\Phi^\pm_\ga$.
\begin{thm} \label{thm:properties Phi}
For $\ga \in \mathcal S_{\mathrm{reg}}$ the functions $\Phi_\ga^\pm$ satisfy the following properties:

\textnormal{(a)} $\Phi^\pm_\ga \in V_{\mu(\ga)}$. 

\textnormal{(b)} For $|\ga|<1$ we have
\[
\int_{\infty(z_-)}^0 \left|\Phi_\ga^-(x)\right|^2 w(x) d_qx<\infty, \qquad 
\int_0^{\infty(z_+)} \left|\Phi_\ga^+(x)\right|^2 w(x) d_qx<\infty.
\]

\textnormal{(c)} The Casorati determinant $v(\ga)=D(\Phi^+_\ga, \Phi_\ga^-)$ is constant on $\R_q$, and
\[
\begin{split}
&v(\ga) = v(\ga;a,b,c,d;z_-,z_+|q) \\
&=  -\frac{z_+(1-q) \te(z_-/z_+) }{\te(cz_-, dz_-, cz_+, dz_+)} 
 \frac{ (cq\ga/as, dq\ga/as, cq\ga/bs, dq\ga/bs, s\ga, q\ga/s;q)_\infty  \te(absz_- z_+/q\ga) }{\ga(q\ga^2;q)_\infty^2}
\end{split}
\]
\end{thm}
\begin{proof}
Properties (a) and (b) follow directly from Proposition \ref{prop:Phi+-} and the asymptotic behavior of $\Phi_\ga^\pm(x)$ for $|x| \rightarrow \infty$, so we only need to check the third property. 

Let $\ga \in \mathcal S_{\mathrm{reg}}\setminus \{sq^k \ | \ k \in \Z_{\geq 0} \}$. Since $\Phi^\pm_\ga \in V_{\mu(\ga)}$ the Casorati determinant $D(\Phi^-_\ga, \Phi_\ga^+)$ is constant on $\R_q$ by Lemma \ref{lem:dim V}. To calculate the determinant we use Proposition \ref{prop:Phi+-}, then
\[
D(\Phi^+_\ga, \Phi^-_\ga) = d_{z_-}(\ga) D(\Phi^+_\ga,\varphi_\ga) + d^\dagger_{z_-}(\ga)D(\Phi^+_\ga,\varphi^\dagger_\ga).
\]
We find from Proposition \ref{prop:c function exp} and Lemma \ref{lem:det Phi},
\[
\begin{split}
D(\Phi_\ga^+,\varphi_\ga) &= c_{z_+}(1/\ga) D(\Phi_\ga^+,\Phi_{1/\ga}^+) = (\ga-1/\ga)c_{z_+}(1/\ga) K_{z_+}, \\
D(\Phi_\ga^+,\varphi_\ga^\dagger) &= c_{z_+}^\dagger(1/\ga) D(\Phi_\ga^+,\Phi_{1/\ga}^+) = (\ga-1/\ga)c_{z_+}^\dagger(1/\ga) K_{z_+},
\end{split}
\]
so we have
\[
D(\Phi^+_\ga, \Phi^-_\ga) = (\ga-1/\ga)K_{z_+} \Big(d_{z_-}(\ga)c_{z_+}(1/\ga)  + d^\dagger_{z_-}(\ga)c_{z_+}^\dagger(1/\ga) \Big) .
\]
From the explicit expression for $d_{z_-}(\ga)$ and $c_{z_+}(\ga)$ we obtain
\[
\begin{split}
d_{z_-}(\ga)c_{z_+}(1/\ga) & + d^\dagger_{z_-}(\ga)c_{z_+}^\dagger(1/\ga) = 
\frac{bsz_- (cq\ga/as, dq\ga/as, cq\ga/bs, dq\ga/bs, s\ga;q)_\infty }{q\ga (\ga^2,q\ga^2,s/\ga;q)_\infty \te(cz_-,dz_-,b/a,az_+,bz_+)  } \\
& \times \Big[ \te(bz_-, az_+, asz_-/\ga, bsz_+/\ga) - \te(az_-, bz_+, bsz_-/\ga, asz_+/\ga) \Big].
\end{split}
\]
Using the $\te$-product identity \eqref{eq:theta prod} with
\begin{align*}
x&=\frac{ise^{i(\al+\be)/2} \sqrt{ |abz_+z_-|}}{\ga}, & y &= ie^{i(\al+\be)/2}\sqrt{|abz_+z_-|}, \\
v&= ie^{i(\al-\be)/2}\sqrt{\left|\frac{ az_-}{bz_+}\right|}, &
w&= ie^{i(\be-\al)/2}\sqrt{\left|\frac{bz_-}{az_+}\right|}.
\end{align*}
where $a=|a|e^{i\al}$ and $b=|b|e^{i\be}$, the term between square bracket equals
\[
bz_+ \te(z_-/z_+, a/b, s/\ga, absz_-z_+/\ga).
\]
Using the explicit expression for $K_{z_+}$ we now find the Casorati determinant given in the theorem. By continuity in $\ga$, the result holds for all $\ga \in \mathcal S_{\mathrm{reg}}$.
\end{proof}

\section{The spectral measure} \label{sec4}
In this section we calculate explicitly the spectral measure $E$ for the self-adjoint operator $(L,\mathcal D)$ using the formula, see \cite[Thm.XII.2.10]{DS63},
\begin{equation} \label{eq:Stieltjes-Perron}
\langle E(\la_1,\la_2)f, g \rangle_{\mathcal L^2} = \lim_{\de \downarrow 0} \lim_{\eps \downarrow 0} \frac{1}{2\pi i} \int_{\la_1 + \de}^{\la_2-\de} \Big( \langle R(\mu+i\eps) f,g \rangle_{\mathcal L^2} - \langle R(\mu-i\eps) f,g \rangle_{\mathcal L^2} \Big) d\mu,
\end{equation}
for $\la_1<\la_2$ and $f,g \in \mathcal L^2$. Here $R(\mu)=(L-\mu)^{-1}$, $\mu \in \C \setminus \R$, denotes the resolvent operator. Our first goal is to find a useful description for the resolvent $R(\mu)$.

\subsection{The resolvent}
Let $\la \in \C \setminus \R$, and let $\ga_\la$ be the unique complex number such that $|\ga_\la|<1$ and $\la=\mu(\ga_\la)$. Note that $\ga_\la \not\in \R$, so $\ga_\la \in \mathcal S_{\mathrm{reg}}$. Let $\mathcal V$ denote the set of zeros of $v(\ga)$, i.e.,
\[
\mathcal V = \left(\bigcup_{\al \in \{\frac{cq}{as}, \frac{dq}{as}, \frac{cq}{bs}, \frac{dq}{bs}, s, \frac{q}{s} \} } \Big\{ \frac{1}{\al q^k}\ \big| \ k \in \Z_{\geq 0} \Big\} \cup \Big\{ z_- z_+q^k\sqrt{abcd/q} \ \big| \ k \in \Z \Big\} \right).
\] 
If $v(\ga)=D(\Phi_\ga^+,\Phi_\ga^-)\neq 0$ the functions $\Phi_\ga^+$ and $\Phi_\ga^-$ are linearly independent, hence for $\ga \in \mathcal S_{\mathrm{reg}} \setminus \mathcal V$ they form a basis for the solution space $V_{\mu(\ga)}$. 

For $\la \in \C\setminus (\R \cup \mu(\mathcal V))$ we define the operator $R_\la:\mathcal D \rightarrow F(\R_q)$ by
\[
(R_\la f)(y) = \langle f, \overline{K_\la(\cdot,y)} \rangle_{\mathcal L^2}, \qquad f\in \mathcal D, \  y \in \R_q, 
\]
where $K_\la:\R_q \times \R_q \rightarrow \C$ is the Green kernel defined by
\[
K_\la(x,y) = 
\begin{cases}
\displaystyle \frac{ \Phi_{\ga_\la}^-(x) \Phi_{\ga_\la}^+(y) }{ v(\ga_\la)}, & x \leq y ,\\ \\
\displaystyle \frac{ \Phi_{\ga_\la}^-(y) \Phi_{\ga_\la}^+(x) }{ v(\ga_\la)}, & x >y.
\end{cases}
\]
Observe that by Theorem \ref{thm:properties Phi} we have $K_\la(x,\cdot) \in \mathcal D$ as well as $K_\la(\cdot,y) \in \mathcal D$ for $x,y \in \R_q$. So $R_\la$ is well-defined as an operator mapping from $\mathcal D$ to $F(\R_q)$. From Propositions \ref{prop:Phi+-} and \ref{prop:Phi pol} we know that the functions $\Phi^\pm_\ga(x)$, considered as functions in $\ga$, are analytic on $\mathcal S_{\mathrm{reg}}$. Now we see that, for $x,y \in \R_q$, the Green kernel $K_{\mu(\ga)}(x,y)$ is a meromorphic function in $\ga$, with poles coming from the zeros of $v(\ga)$. 

\begin{prop}
For $\la \in \C \setminus (\R \cup \mu(\mathcal V))$, the operator $R_{\la}$ is the resolvent of $(L,\mathcal D)$.
\end{prop}
\begin{proof}
The operator $(L,\mathcal D)$ is self-adjoint, hence the spectrum is contained in $\R$. So for $\la \in \C\setminus \R$ the resolvent $R(\la)$ is a bounded linear operator mapping from $\mathcal L^2$ to $\mathcal D$, and therefore for a given $y \in \R_q$ the assignment $f \mapsto (R(\la)f)(y)$ defines a bounded linear functional on $\mathcal L^2$. By the Riesz representation theorem there exists a kernel $K'_\la(\cdot,y) \in \mathcal L^2$ such that $(R(\la) f)(y) = \langle f, K'_\la(\cdot,y) \rangle_{\mathcal L^2}$. So it suffices to show that $(L-\la)R_\la f = f$ for $f \in \mathcal D$.

Suppose that $y> 0$, then, 
\[
\begin{split}
(&L-\la)R_\la f(y)\\
= & \int_{\R_q} f(x) \Big(A(y) K_\la(x,y/q) + B(y) K_\la(x,yq) + (C(y)-\la) K_\la(x,y)\Big) w(x) d_q x \\
= & \frac{1}{v({\ga_\la})}\int_{\infty(z_-)}^{yq} f(x)\Phi_{\ga_\la}^-(x) \Big(A(y)\Phi_{\ga_\la}^+(y/q) + B(y)\Phi_{\ga_\la}^+(yq) + (C(y)-\la)\Phi_{\ga_\la}^+(y) \Big) w(x)d_qx \\
+&\frac{1}{v({\ga_\la})}\int_{y/q}^{\infty(z_+)} f(x)\Phi_{\ga_\la}^+(x) \Big(A(y)\Phi_{\ga_\la}^-(y/q) + B(y)\Phi_{\ga_\la}^-(yq) + (C(y)-\la)\Phi_{\ga_\la}^-(y) \Big) w(x)d_qx \\
+&\frac{(1-q)yw(y)}{v({\ga_\la})} f(y) \Big(A(y)\Phi_{\ga_\la}^-(y) \Phi_{\ga_\la}^+(y/q)+B(y)\Phi_{\ga_\la}^-(yq) \Phi_{\ga_\la}^+(y) + (C(y)-\la)\Phi_{\ga_\la}^-(y) \Phi_{\ga_\la}^+(y)  \Big)\\
=& \frac{(1-q)yB(y)w(y)}{v({\ga_\la})} f(y) \Big(\Phi_{\ga_\la}^+(y)\Phi_{\ga_\la}^-(yq) - \Phi_{\ga_\la}^+(yq)\Phi_{\ga_\la}^-(y) \Big)\\
=& f(y).
\end{split}
\]
Here we used that $\Phi_{\ga_\la}^\pm$ are solutions of $Lf=\la f$, $v(\ga) = D(\Phi_\ga^+,\Phi_\ga^-)$, and Definition \ref{def:Casorati} of the Casorati determinant.  The proof for $y < 0$ runs along the same lines.
\end{proof}

\subsection{The continuous spectrum}
We are going to investigate the integrand in \eqref{eq:Stieltjes-Perron}. Using the definition of the Green kernel we have
\begin{equation} \label{eq:int R}
\begin{split}
\langle R_\mu f, g \rangle_{\mathcal L^2} &= \iint \limits_{\R_q \times \R_q} f(x)\overline{g(y)} K_\mu(x,y) w(x)w(y) d_qx\, d_qy\\
&= \iint \limits_{\substack{(x,y) \in  \R_q \times \R_q \\ x \leq y }} \frac{ \Phi^-_{\ga_\mu}(x) \Phi^+_{\ga_\mu}(y)}{v(\ga_\mu)} \big(f(x)\overline{g(y)}+f(y)\overline{g(x)}\big)\big(1-\hf\de_{x,y}\big) w(x) w(y) d_qx\,d_qy.
\end{split}
\end{equation}
The Kronecker-delta function $\de_{xy}$ is needed here to prevent the terms on the diagonal $x=y$ from being counted twice. 

We define two functions $v_1$ and $v_2$ that we need to describe the spectral measure $E$;
\begin{align*}
v_1(\ga) &= \frac{ (cq/a, dq/a;q)_\infty^2 \te(bz_+,bz_-) }{(1-q)abz_-^2z_+^2 \te(z_-/z_+,z_+/z_-,a/b,b/a)}\\
& \times \frac{(\ga^{\pm2};q)_\infty}{(s\ga^{\pm1}, cq\ga^{\pm1}/as, dq\ga^{\pm1}/as;q)_\infty \te(s\ga^{\pm1},absz_- z_+\ga^{\pm 1})}\\
& \times \Big( z_-\te(az_+, cz_+, dz_+, bz_-, asz_-\ga^{\pm1}) -  z_+\te(az_-, cz_-, dz_-, bz_+, asz_+\ga^{\pm1})  \Big),\\
v_2(\ga) &= \frac{(cq/a,dq/a,cq/b,dq/b;q)_\infty \te(az_+, az_-, bz_+, bz_-, cdz_- z_+) }{abz_-^2z_+(1-q)\te(z_+/z_-,a/b,b/a)} \\
& \times \frac{ (\ga^{\pm2};q)_\infty }{ (s\ga^{\pm1};q)_\infty \te(s\ga^{\pm1}, absz_- z_+\ga^{\pm1}) }.
\end{align*}
Note that $v_1$ and $v_2$ are both invariant under $\ga \leftrightarrow 1/\ga$. Let $\mathcal D_{\text{fin}} \subset \mathcal D$ be the subspace consisting of finitely supported functions in $\mathcal L^2$. To a function $f \in \mathcal D_{\mathrm{fin}}$ we associate two functions $\mathcal F_cf$ and $\mathcal F_c^\dagger f$ on the unit circle $\T=\{z\in \C \ | \ |z|=1\}$ defined by
\[
\begin{split}
(\mathcal F_cf)(\ga) &= \int_{\R_q} f(x) \varphi_\ga(x) w(x) d_q x, \\
(\mathcal F_c^\dagger f)(\ga) &= \int_{\R_q} f(x) \varphi_\ga^\dagger(x) w(x) d_q x,
\end{split}
\]
where $\ga \in \T$. 

We are now almost ready to describe the spectral measure $E((\la_1,\la_2))$ for $(\la_1,\la_2) \subset (-2,2)$. First we give a preliminary result. The proof is an easy, but rather tedious computation that we  carry out in the appendix.
\begin{lem} \label{lem:I(x,y)}
For $x,y \in \R_q$ and  $\ga, \ga^{-1} \in \mathcal S_{\mathrm{reg}}\setminus \mathcal V$ we have
\[
\begin{split}
&\frac{ \Phi_{1/\ga}^-(x) \Phi_{1/\ga}^+(y)}{v(1/\ga)} - \frac{ \Phi_{\ga}^-(x) \Phi_{\ga}^+(y)}{v(\ga)} = \\ 
&\quad\frac{1}{\ga-1/\ga}\Big[ v_1(\ga) \varphi_\ga(x) \varphi_\ga(y) + v_2(\ga) \big(\varphi_\ga(x) \varphi_\ga^\dagger(y) + \varphi_\ga^\dagger(x) \varphi_\ga(y) \big) + v_1^\dagger(\ga) \varphi_\ga^\dagger(x)\varphi_\ga^\dagger(y)\Big].
\end{split}
\]
\end{lem}

\begin{prop} \label{prop:E(-2,2)}
Let $(a,b,c,d) \in P_{\mathrm{gen}}$, let $0< \psi_1<\psi_2< \pi$, and let $\la_1=\mu(e^{i\psi_2})$ and $\la_2=\mu(e^{i\psi_1})$. Then for $f,g \in \mathcal D_{\text{fin}}$,
\[
\begin{split}
\langle E(\la_1,\la_2)f,g \rangle_{\mathcal L^2}= \frac{1}{2\pi}\int_{\psi_1}^{\psi_2}\big( (\mathcal F_c^\dagger \overline{g})(e^{i\psi}) \, (\mathcal F_c \overline{g})(e^{i\psi}) \big)
\begin{pmatrix}
v_2(e^{i\psi}) & v_1^\dagger(e^{i\psi}) \\
v_1(e^{i\psi}) & v_2(e^{i\psi})
\end{pmatrix}
\vect{(\mathcal F_cf)(e^{i\psi})}{(\mathcal F_c^\dagger f)(e^{i\psi})} d\psi.
\end{split}
\]
\end{prop}
\begin{proof}
Let $\la \in (-2,2)$, then $\la = \mu(e^{i\psi})$ for a unique $\psi \in (0,\pi)$.  In this case we have
\[
\lim_{\eps \downarrow 0} \ga_{\la \pm i\eps} = e^{\mp i\psi}.
\]
Now we obtain
\[
\begin{split}
\lim_{\eps \downarrow 0} &\Big( \frac{ \Phi_{\ga_{\la+i\eps}}^-(x) \Phi_{\ga_{\la+i\eps}}^+(y)}{v(\ga_{\la+i\eps})} - \frac{ \Phi_{\ga_{\la-i\eps}}^-(x) \Phi_{\ga_{\la-i\eps}}^+(y)}{v(\ga_{\la-i\eps})} \Big) \\
&= \frac{ \Phi_{e^{-i\psi}}^-(x) \Phi_{e^{-i\psi}}^+(y)}{v(e^{-i\psi})} - \frac{ \Phi_{e^{i\psi}}^-(x) \Phi_{e^{i\psi}}^+(y)}{v(e^{i\psi})},
\end{split}
\]
which is symmetric in $x$ and $y$ by Lemma \ref{lem:I(x,y)}. Symmetrizing the double $q$-integral from \eqref{eq:int R} then gives
\[
\begin{split}
\lim_{\eps \downarrow 0}&\Big(\langle R_{\la+i\eps} f,g \rangle_{\mathcal L^2} - \langle R_{\la-i\eps} f,g \rangle_{\mathcal L^2}\Big) = \\
&\iint \limits_{\R_q \times \R_q}  
\frac{1}{\ga-1/\ga}\Big[ v_1(\ga) \varphi_\ga(x) \varphi_\ga(y) + v_2(\ga) \big(\varphi_\ga(x) \varphi_\ga^\dagger(y) + \varphi_\ga^\dagger(x) \varphi_\ga(y) \big) + v_1^\dagger(\ga) \varphi_\ga^\dagger(x)\varphi_\ga^\dagger(y)\Big] \\
& \quad \times  f(x) \overline{g(y)} w(x) w(y) d_qx\, d_qy,
\end{split}
\]
where $\ga = e^{i\psi}$.
Rewriting this expression in vector notation and using formula \eqref{eq:Stieltjes-Perron}, we obtain the desired result.
\end{proof}
The previous proposition implies that $(-2,2)$ is contained in the spectrum $\si(L)$ of $L$. Since $\varphi_\ga, \varphi_\ga^\dagger \not\in \mathcal L^2$ for $\ga \in \T$, $(-2,2)$ is part of the continuous spectrum. Observe that the spectral projection is on a 2-dimensional space of eigenvectors, so $(-2,2)$ has multiplicity two. Because the spectrum is a closed set, the points $-2$ and $2$ must be elements of the spectrum $\si(L)$.
\begin{lem} \label{lem:-2 2}
The points $-2$ and $2$ are elements of the continuous spectrum of $L$.
\end{lem}
\begin{proof}
Since the residual spectrum of a self-adjoint operator is empty, $\mu(-1)=-2$ and $\mu(1)=2$ must either be elements of the point spectrum or the continuous spectrum. We show that $2$ is not in the point spectrum of $L$. The proof for $-2$ is the same.

Suppose that there exists a function $f \in \mathcal L^2$ that satisfies $Lf=2f$, then the restriction $f^{\mathrm{res}}$ of $f$ to $\R_q^+$ is an element of $V_2^+$. From Lemma \ref{lem:ga = 1} it follows that $f^{\mathrm{res}}=\al \Phi_1^+ + \be \left.\frac{ d\Phi_\ga^+}{d\ga} \right|_{\ga=1}$ for some coefficients $\al$ and $\be$. But neither of the functions $\Phi_1^+$ and $\left.\frac{ d\Phi_\ga^+}{d\ga} \right|_{\ga=1}$ is integrable with respect to $w(x)$ on $\R_q^+$, see \eqref{eq:asymp Phi} and Lemma \ref{lem:asymp w&u}, which contradicts the fact that $f \in \mathcal L^2$.
\end{proof}

\subsection{The point spectrum}
Let $\mu \in \R\setminus [-2,2]$, then 
\[
\lim_{\eps \downarrow 0} \ga_{\mu\pm i \eps} = \ga_\mu.
\]
From \eqref{eq:Stieltjes-Perron} and \eqref{eq:int R} we see that in this case the only contribution to the spectral measure $E$ comes from the real poles of the Green kernel $K_{\mu(\ga)}(x,y)$, $x,y \in \R_q$, considered as a function of $\ga$. Let $\Ga \subset \mathcal V$ denote the set of poles of the Green kernel inside the interval $(-1,1)$. We now have the following property for the spectral measure.
\begin{prop} \label{prop:E=0}
For real numbers $\mu_1< \mu_2$ satisfying \mbox{$(\mu_1,\mu_2) \cap \Big( \mu(\Ga) \cup [-2,2] \Big) = \emptyset$}, we have $E((\mu_1,\mu_2))=0$.
\end{prop}
The set $\Ga$ of real poles of the Green kernel inside the unit disc is given by
\[
\begin{split}
\Ga &=  \Ga^{\mathrm{fin}}_{s} \cup \Ga^{\mathrm{fin}}_{q/s} \cup \Ga^{\mathrm{fin}}_{dq/as}\cup \Ga^{\mathrm{inf}},\\
\Ga^{\mathrm{fin}}_\al &=  \Bigg\{ \frac{1}{\al q^{k} }\ \Big| \ k\in\Z_{\geq 0},\ \al q^{k}>1 \Bigg\},\\
\Ga^{\mathrm{inf}} &= \Bigg\{ z_-z_+q^{k}\sqrt{abcd/q} \ \Big| \ k \in \Z, \ -z_-z_+q^k\sqrt{abcd/q}<1 \Bigg\}.
\end{split}
\]
The superscripts `fin' and `inf' refer to the finite or infinite cardinality of the sets. Recall that for $a,b,c,d \in \R$ we have assumed that $q<a/b<1$ and $q<c/d<1$, therefore
\[
 \frac{acq}{bd} ,  \frac{bcq}{ad} ,  \frac{adq}{bc} <1, \qquad 1<\frac{bdq}{ac}<\frac{1}{q}.
\]
So the factor $(cq\ga/as, dq\ga/as, cq\ga/bs, dq\ga/bs;q)_\infty$ of the function $v(\ga)$ has at most one zero inside the interval $(-1,1)$, and this zero only occurs when $a,b,c,d \in \R$. This shows that the set $\Ga_{dq/as}^{\mathrm{fin}}$ has at most one element. We remark that for $(a,b,c,d) \in P_{\mathrm{gen}}$ the real poles of the Green kernel are simple. 

Next we need to find the spectral measure on the set $\mu(\Ga)$. For this the following lemma is useful.
\begin{lem} \label{lem:Phi on si_p}
For $\ga \in \Ga$ we have 
\[
\Phi_\ga^+(x) = b(\ga) \Phi_\ga^-(x), \qquad x \in \R_q,
\]
where $b(\ga)=b(\ga;a,b,c,d;z_-,z_+;q)$ is given by
\[
b(\ga)=
\begin{cases}
\left(\dfrac{z_+}{z_-}\right)^{k+1}\dfrac{\te(cz_-,dz_-)}{\te(cz_+,dz_+)}, &\ga =z_-z_+q^{k}\sqrt{abcd/q} \in \Ga^{\mathrm{inf}},\\
\left(\dfrac{z_-}{z_+}\right)^k ,& \ga = s^{-1}q^{-k} \in \Ga^{\mathrm{fin}}_s,\\
\left( \dfrac{z_-}{z_+} \right)^k \dfrac{ \te(az_+, bz_+, cz_-, dz_-) }{ \te(az_-, bz_-, cz_+, dz_+) }, &\ga = sq^{-1-k} \in \Ga^{\mathrm{fin}}_{q/s},\\
\dfrac{ \te(bz_+,cz_-)}{\te(bz_-,cz_+)}, &\ga = \dfrac{as}{dq} \in \Ga^{\mathrm{fin}}_{dq/as}.
\end{cases}
\]
\end{lem}
\begin{proof}
If $\ga \in \Ga$, then $v(\ga) = D(\Phi_\ga^-,\Phi_\ga^+) = 0$, hence $\Phi_\ga^+ = b(\ga) \Phi_\ga^-$ for some nonzero factor $b(\ga)$. For $\ga_k = s^{-1}q^{-k} \in \Ga^{\mathrm{fin}}_s$ the value of $b(\ga_k)$ follows from Corollary \ref{cor:Phi pol}. For the other cases it is enough by Proposition \ref{prop:Phi+-} to show that $d_{z_+}(\ga)=b(\ga) d_{z_-}(\ga)$ and $d_{z_+}^\dagger(\ga)=b(\ga) d_{z_-}^\dagger(\ga)$. This is verified by a straightforward calculation. Note that for $\ga = as/dq \in \Ga^{\mathrm{fin}}_{dq/as}$ we have $d_{z_+}^\dagger(\ga)=d_{z_-}^\dagger(\ga)=0$.
\end{proof}

We are now ready to calculate the spectral measure $E$ of $(L,\mathcal D)$ for the discrete part of the spectrum $\si(L)$. We will write $E(\{\mu(\ga)\})$ for the spectral measure $E((a,b))$, if $(a,b)$ is an interval such that $(a,b) \cap \mu(\Ga) = \{\mu(\ga)\}$. For $f \in \mathcal L^2$ we define a function $\mathcal F_pf$ on $\Ga$ by
\[
(\mathcal F_pf)(\ga) = \langle f, \Phi_\ga^+ \rangle_{\mathcal L^2}, \qquad \ga \in \Ga.
\]
Note that Theorem \ref{thm:properties Phi}(b) and Lemma \ref{lem:Phi on si_p} imply that $\Phi_\ga^+ \in \mathcal L^2$ for $\ga \in \Ga$, so the inner product above exists for all $f \in \mathcal L^2$.
\begin{prop} \label{prop:E point}
Let $(a,b,c,d) \in P_{\mathrm{gen}}$. For $f,g \in \mathcal L^2$ and $\ga \in \Ga$, the spectral measure $E(\{\mu(\ga)\})$ is given by
\[
\langle E(\{\mu(\ga)\})f,g \rangle_{\mathcal L^2} = (\mathcal F_pf)(\ga) \overline{(\mathcal F_p g)(\ga)} N(\ga),
\]
where
\[
N(\ga)=N(\ga;a,b,c,d;z_-,z_+|q) = b(\ga)^{-1}\Res{\la=\ga} \left(\frac{ 1/\la-\la}{\la v(\la)} \right),
\]
and $b(\ga)$ is given in Lemma \ref{lem:Phi on si_p}. 
\end{prop}
\begin{proof}
Let $\ga \in \Ga$, and $f,g \in \mathcal L^2$. We use \eqref{eq:Stieltjes-Perron} and \eqref{eq:int R} to calculate the spectral measure $E(\{\mu(\ga)\})$. By the residue theorem we find
\[
\begin{split}
\langle E(\{\mu(\ga)\})f,g \rangle_{\mathcal L^2} =& \frac{1}{2\pi i}\int_{\mathcal C} \langle R_\mu f,g \rangle_{\mathcal L^2} d\mu \\
=&\iint \limits_{\substack{(x,y) \in  \R_q \times \R_q \\ x \leq y }} -\Res{\la= \ga} \left((1-1/\la^2)\frac{ \Phi^-_{\la}(x) \Phi^+_{\la}(y)}{v(\la)}\right) \\ 
&\qquad \times \big(f(x)\overline{g(y)}+f(y)\overline{g(x)}\big)\big(1-\hf\de_{x,y}\big) w(x) w(y) d_qx\,d_qy
\end{split}
\]
Here $\mathcal C$ is a clockwise oriented contour encircling $\mu(\ga)$ once, and $\mathcal C$ does not encircle any other points in $\Ga$. The factor $1-1/\la^2$ comes from changing the integration variable $\mu=\mu(\la)$ to $\la$. By Lemma \ref{lem:Phi on si_p} we have $\Phi_{\ga}^+ = b(\ga) \Phi_{\ga}^-$, so we may symmetrize the double $q$-integral, and then
\[
\langle E(\{\mu(\ga)\})f,g \rangle_{\mathcal L^2} = b(\ga)^{-1} \langle f,\Phi_{\ga}^+\rangle_{\mathcal L^2} \langle \Phi_{\ga}^+, g \rangle_{\mathcal L^2} \, \Res{\la=\ga} \left(\frac{ 1/\la-\la}{\la v(\la)} \right). 
\]
This proves the proposition.
\end{proof}
It is an easy exercise to calculate the weight $N(\ga)$, $\ga \in \Ga$, explicitly. The result is as follows. For $\ga= as/dq \in \Ga^{\mathrm{fin}}_{dq/as}$ we have
\[
N(\ga) = \frac{ \te(bz_-, cz_+, cz_+, dz_-, dz_+) (ac/bdq;q)_\infty}{ z_+(1-q) \te(bz_+, bdz_-z_+, z_-/z_+) (q,a/b, a/d, c/b, c/d;q)_\infty }, 
\]
for $\ga = s^{-1}q^{-k} \in \Ga^{\mathrm{fin}}_{s}$
\[
\begin{split}
N(\ga) =& \frac{ \te(cz_-,cz_+, dz_-,dz_+) (ab/cdq;q)_\infty }{ z_+(1-q)\te(z_-/z_+, cdz_-z_+) (q,a/c,a/d, b/c,b/d;q)_\infty } \\
& \times \frac{ (q^2cd/ab, qcd/ab;q)_{2k} }{ (q,cq/a, cq/b, dq/a, dq/b, qcd/ab;q)_k } \left(- \frac{abz_+^2}{q} \right)^k q^{\frac32 k(k-1)},
\end{split}
\]
for $\ga = sq^{-1-k} \in \Ga^{\mathrm{fin}}_{q/s}$
\[
\begin{split}
N(\ga) = & \frac{ \te(az_-,bz_-, cz_+, cz_+, dz_+, dz_+) (cd/abq;q)_\infty }{z_+(1-q) \te(az_+, bz_+, z_-/z_+, abz_-z_+) (q, c/a,c/b, d/a,d/b;q)_\infty } \\
& \times \frac{ (qab/cd, ab/cd;q)_{2k} }{ (q,aq/c,aq/d, bq/d, bq/d, abq/cd;q)_k } \left( -\frac{ z_+^2 q^4 }{cd } \right)^k q^{\frac32 k(k-1)},
\end{split}
\]
and finally for $\ga = absz_-z_+q^{k-1} \in \Ga^{\mathrm{inf}}$ we have
\[
\begin{split}
N(\ga) = & \frac{ \te(cz_+, dz_+)^2 (abcdz_-^2z_+^2/q, abcdz_-^2z_+^2;q)_\infty }{z_+(1-q) \te(z_-/z_+) (q,q,abz_-z_+,acz_-z_+, adz_-z_+, bcz_-z_+, bdz_-z_+, cdz_-z_+;q)_\infty } \\
& \times \frac{ (abz_-z_+,acz_-z_+, adz_-z_+, bcz_-z_+, bdz_-z_+, cdz_-z_+;q)_k }{ (abcdz_-^2z_+^2/q, abcdz_-^2z_+^2;q)_{2k} } \left( \frac{z_-}{z_+} \right)^{k+1} (-1)^kq^{\hf k(k+1)}.
\end{split}
\]

As a result of Proposition \ref{prop:E point} we obtain orthogonality relations for $\Phi_\ga^+$, $\ga \in \Ga$.
\begin{cor} \label{cor:orth rel}
Let $\ga, \ga' \in \Ga$, then
\[
\langle \Phi_\ga^+, \Phi_{\ga'}^+ \rangle_{\mathcal L^2} = \frac{\de_{\ga \ga'}}{N(\ga)}.
\]
\end{cor}
\begin{proof}
Eigenfunctions corresponding to different eigenvalues of a self-adjoint operator are pairwise orthogonal. Since for $\ga, \ga' \in \Ga$, $\ga \neq \ga'$, the functions $\Phi_\ga^+$ and $\Phi_{\ga'}^+$ are eigenfunctions of $(L,\mathcal D)$ with distinct eigenvalues $\mu(\ga)$ and $\mu(\ga')$, orthogonality follows.

Let $\ga \in \Ga$. By Proposition \ref{prop:E point},
\[
\langle \Phi_\ga^+, \Phi_{\ga}^+ \rangle_{\mathcal L^2} = \langle E(\{\mu(\ga)\})\Phi_\ga^+, \Phi_{\ga}^+ \rangle_{\mathcal L^2} =N(\ga) \langle \Phi_\ga^+, \Phi_{\ga}^+ \rangle_{\mathcal L^2}^2 ,
\]
from which the squared norm of $\Phi_\ga^+$ follows.
\end{proof}
\begin{rem}
For $\ga, \ga' \in \Ga^{\mathrm{fin}}_s$ Corollary \ref{cor:orth rel} gives orthogonality relations for a finite number of big $q$-Jacobi polynomials, see Proposition \ref{prop:Phi pol} and Lemma \ref{lem:phi pol}.
\end{rem}
Since $\Phi_{s^{-1}}^+(x)=1$, Corollary \ref{cor:orth rel} gives an evaluation of the integral $\langle 1,1 \rangle_{\mathcal L^2}$ in case $\Ga^{\mathrm{fin}}_s$ is not empty, i.e., if $s>1$. 
\begin{cor}
For $\sqrt{ab/cdq}<1$ we have
\[
\begin{split}
\frac{1}{1-q}&\int_{\R_q} \frac{ (ax, bx;q)_\infty }{ (cx,dx;q)_\infty } d_q x \\
=& z_+ \frac{ (az_+, bz_+;q)_\infty }{ (cz_+, dz_+;q)_\infty } {}_2\psi_2 \left( \genfrac{.}{.}{0pt}{}{cz_+,dz_+}{az_+,bz_+}\ ;q,q \right)- z_- \frac{ (az_-, bz_-;q)_\infty }{ (cz_-, dz_-;q)_\infty } {}_2\psi_2 \left( \genfrac{.}{.}{0pt}{}{cz_-,dz_-}{az_-,bz_-}\ ;q,q \right) \\
=& z_+\frac{  (q,a/c,a/d,b/c,b/d;q)_\infty \te(z_-/z_+, cdz_-z_+) }{ (ab/cdq;q)_\infty \te(cz_-, dz_-, cz_+, dz_+;q)_\infty },
\end{split}
\]
\end{cor}
Here $_2\psi_2$ denotes the usual bilateral series as defined in \cite{GR04}.
\begin{rem}
This is the summation formula from \cite[Exer.5.10]{GR04}, and it is actually valid without the restrictions on $a,b,c,d$ as long as the denominator of the integrand is nonzero for all $x \in \R_q$. Note that there is a misprint in \cite[Exer.5.10]{GR04}: the factors $(e/ab,q^2f/e;q)_\infty$ on the left hand side must be replaced by $(c/qf,q^2f/c;q)_\infty$.
\end{rem}

\begin{cor}
Let $(a,b,c,d) \in P_{\mathrm{gen}}$. The spectrum of the self-adjoint operator $(L,\mathcal D)$ consists of the continuous spectrum $\si_c(L) = [-2,2]$, with multiplicity two, and the point spectrum $\si_p(L)=\mu(\Ga)$, with multiplicity one. 
\end{cor}
\begin{proof}
This follows from Propositions \ref{prop:E(-2,2)}, \ref{prop:E=0}, \ref{prop:E point} and Lemma \ref{lem:-2 2}.
\end{proof}

\section{The vector-valued big $q$-Jacobi function transform} \label{sec5}
In this section we define the vector-valued big $q$-Jacobi function transform $\mathcal F$, which is closely related to the maps $\mathcal F_c$, $\mathcal F_c^\dagger$ and $\mathcal F_p$. We show that $\mathcal F$ is an isometric isomorphism mapping from $\mathcal L^2$ into a certain Hilbert space $\mathcal H$, and we also determine $\mathcal F^{-1}$. The vector-valued big $q$-Jacobi function transform diagonalizes the second order difference operator $L$; let $M$ be the multiplication operator defined by $(M f)(\ga) = \mu(\ga)f(\ga)$ for all $\mu(\ga) \in \si(L)$, then
\[
(\mathcal F \circ L \circ \mathcal F^{-1})f = M f,
\]
for all $f \in \mathcal H$ such that $M f \in \mathcal H$.

We still assume that $(a,b,c,d) \in P_{\mathrm{gen}}$, and we distinguish between the cases $a \neq \overline{b}$ and $a=\overline{b}$. For a vector $y=\vect{y_1}{y_2} \in \C^2$ we denote
\[
y^T=
\begin{cases}
(\overline{y_1}\ \overline{y_2} ), &\text{if }a = \overline{b},\\
(\overline{y_2}\ \overline{y_1} ), &\text{if }a \neq \overline{b}.
\end{cases}
\]
With this convention we have for $\ga \in \R \cup \T$ and $x \in \R_q$
\[
\vect{\varphi_\ga(x)}{\varphi_\ga^\dagger(x)}^T =  \big( \varphi_\ga^\dagger (x) \ \varphi_\ga(x) \big),
\]
since $\overline{\varphi_\ga(x)} = \varphi_\ga^\dagger(x)$ if $a = \overline{b}$, and $\overline{\varphi_\ga(x)}=\varphi_\ga(x)$ if $a \neq \overline{b}$. 

\subsection{The vector-valued big $q$-Jacobi function transform $\mathcal F$}
Let $F(\T \cup \Ga)$ be the linear space consisting of functions that are complex-valued on $\Ga$ and $\C^2$-valued on $\T$. With the maps $\mathcal F_c$, $\mathcal F_c^\dagger$ and $\mathcal F_p$ we define an integral transform $\mathcal F:\mathcal D_{\mathrm{fin}} \rightarrow F(\T\cup \Ga)$.
\begin{Def}
For $f \in \mathcal D_{\mathrm{fin}}$ we define the vector-valued big $q$-Jacobi function $\mathcal F$ by
\[
(\mathcal F f)(\ga) =
\begin{cases}
\displaystyle \vect{(\mathcal F_c f)(\ga)}{(\mathcal F_c^\dagger f)(\ga)}, & \ga \in \T,\\ \\
(\mathcal F_p f)(\ga), & \ga \in \Ga.
\end{cases}
\]
\end{Def}
We define a kernel $\Psi(x,\ga)$, $x \in \R_q$, $\ga \in \T \cup \Ga$, by
\begin{equation} \label{eq:Psi}
\Psi(x,\ga) =
\begin{cases}
\displaystyle \vect{\varphi_\ga(x)}{\varphi_\ga^\dagger(x)}, & \ga \in \T,\\ \\
\Phi_\ga^+(x), & \ga \in \Ga,
\end{cases}
\end{equation}
We may write $\mathcal F$ as an integral transform with kernel $\Psi$,
\[
(\mathcal Ff)(\ga) = \int_{\R_q} f(x) \Psi(x,\ga) w(x) d_qx, \qquad f \in \mathcal D_{\mathrm{fin}},\quad \ga \in \T \cup \Ga.
\]
We are going to show that $\mathcal F$ extends to a continuous operator mapping from $\mathcal L^2$ into a Hilbert space $\mathcal H$, that we now define. 

We define a matrix-valued function $\mathbf v$ on $\T$ by
\[
\ga \mapsto \mathbf v(\ga) = 
\begin{pmatrix}
v_2(\ga) & v_1^\dagger(\ga) \\
v_1(\ga) & v_2(\ga)
\end{pmatrix}.
\]
We remark that $\mathbf v(\ga)$, $\ga \in \T$, is positive-definite. Let $\mathcal H_c'$ be the Hilbert space consisting of $\C^2$-valued functions on $\T$, that have finite norm with respect to the inner product
\[
\langle g_1, g_2 \rangle_{\mathcal H_c'} = \frac{1}{4\pi i} \int_\T g_2(\ga)^T \mathbf v(\ga)
g_1(\ga) \frac{ d\ga}{\ga},
\]
where the unit circle $\T$ is oriented in the counter-clockwise direction. Let $r$ denote the reflection operator defined by $(r g)(\ga) = g(\ga^{-1})$. We define the Hilbert space $\mathcal H_c$ to be the subspace of $\mathcal H_c'$ consisting of functions $g$ that satisfy $r g=g$ in $\mathcal H_c'$. We denote the inner product on $\mathcal H_c$ by $\langle \cdot, \cdot \rangle_{\mathcal H_c}$. Furthermore, let $\mathcal H_p$ be the Hilbert space consisting of complex-valued functions on $\Ga$, that have finite norm with respect to the inner product 
\[
\langle g_1, g_2 \rangle_{\mathcal H_p} = \sum_{\ga \in \Ga} g_1(\ga) \overline{g_2(\ga)} N(\ga).
\]
We define the Hilbert space $\mathcal H \subset F(\T\cup \Ga)$ by $\mathcal H = \mathcal H_c \oplus \mathcal H_p$. 
\begin{prop} \label{prop:F iso}
The map $\mathcal F$ extends uniquely to an operator $\mathcal F : \mathcal L^2 \rightarrow \mathcal H$, satisfying
\[
\langle \mathcal Ff_1, \mathcal Ff_2 \rangle_{\mathcal H} = \langle f_1,f_2 \rangle_{\mathcal L^2}, \qquad f_1,f_2 \in \mathcal L^2.
\]
Hence, $\mathcal F$ is an isometric isomorphism onto its range $\mathcal R(\mathcal F) \subset \mathcal H$.
\end{prop}
\begin{proof}
Let $f_1,f_2 \in \mathcal D_{\mathrm{fin}}$. Combining Propositions \ref{prop:E(-2,2)} and \ref{prop:E point} we find
\[
\langle f_1, f_2 \rangle_{\mathcal L^2} = \langle E(\R) f_1, f_2 \rangle_{\mathcal L^2} = \langle \mathcal Ff_1, \mathcal Ff_2 \rangle_{\mathcal H}.
\]
Here we used
\[
\big( (\mathcal F_c^\dagger \overline{f})(\ga)\  (\mathcal F_c \overline{f})(\ga) \big) =
\vect{(\mathcal F_c^\dagger f)(\ga)}{(\mathcal F_c f)(\ga)}^T.
\]
Since $\mathcal D_{\mathrm{fin}}$ is dense in $\mathcal L^2$, the map $\mathcal F$ extends uniquely to a continuous operator, also denoted by $\mathcal F$, mapping isometrically into $\mathcal R(\mathcal F) \subset \mathcal H$. 
\end{proof}

\begin{lem} \label{lem:Psi in H}
Let $y \in \R_q$ and let $f_y(x)=\de_{xy}/w(y) \in \mathcal L^2$, then $\mathcal Ff_y = \Psi(y,\cdot) \in \mathcal H$.
\end{lem}
\begin{proof}
We have
\[
(\mathcal F f_y)(\ga) = 
\begin{cases}
\displaystyle \vect{\varphi_\ga(y)}{\varphi_\ga^\dagger(y)}, & \ga \in \T,\\ \\
\Phi^+_\ga(y), & \ga \in \Ga,
\end{cases}
\]
so $\mathcal Ff_y=\Psi(y, \cdot)$. By Proposition \ref{prop:F iso} this lies in $\mathcal H$.
\end{proof}

We define an integral transform $\mathcal G : \mathcal H \rightarrow F(\R_q)$ by
\[
(\mathcal Gg)(x) = \langle g, \Psi(x,\cdot) \rangle_{\mathcal H}, \qquad g \in \mathcal H, \quad x \in \R_q.
\]
By Lemma \ref{lem:Psi in H} this inner product exists for all $g \in \mathcal H$. We denote by $\mathcal G_c$, respectively $\mathcal G_p$, the integral transform $\mathcal G$ restricted to $\mathcal H_c$, respectively $\mathcal H_p$.
\begin{prop} \label{prop:GF=id}
$\mathcal G \mathcal F = id_{\mathcal L^2}$
\end{prop}
\begin{proof}
Let $f \in \mathcal L^2$ and let $f_y \in \mathcal L^2$ be defined as in Lemma \ref{lem:Psi in H}, then it follows from Proposition \ref{prop:F iso} that
\[
f(y) = \langle f, f_y \rangle_{\mathcal L^2} = \langle \mathcal Ff, \mathcal Ff_y \rangle_{\mathcal H} = 
\big(\mathcal G(\mathcal Ff) \big)(y).\qedhere
\]
\end{proof}
We showed that $\mathcal G$ is a left inverse of $\mathcal F$. Next we are going show that $\mathcal G$ is also a right inverse. We do this for the transforms $\mathcal G_c$ and $\mathcal G_p$ separately. First a preliminary result. We denote by $\langle f, g \rangle_{k;n}$ the limit of the truncated inner product $\lim_{l,m \rightarrow \infty} \langle f,g \rangle_{k,l;m,n}$, provided that this limit exists.
\begin{lem} \label{lem:inprod phi psi}
Let $\ga_1,\ga_2 \in \C^*$ such that $\mu(\ga_1) \neq \mu(\ga_2)$, then for $\phi \in V_{\mu(\ga_1)}$ and $\psi \in V_{\mu(\ga_2)}$,
\[
\langle \phi, \overline{\psi} \rangle_{k;n} = \frac{D(\phi, \psi)(z_+q^{n-1})- D(\phi, \psi)(z_-q^{k-1})}{\mu(\ga_1)-\mu(\ga_2)}.
\]
\end{lem}
\begin{proof}
Functions in $V_\mu$, $\mu \in \C$, are continuously differentiable at the origin, therefore
\[
\lim_{l \rightarrow \infty} D(\phi, \psi)(z_\pm q^{l}) =0.
\]
Since $\phi$ and $\psi$ are eigenfunctions of $L$ for eigenvalue $\mu(\ga_1)$, respectively $\mu(\ga_2)$, we obtain from Proposition \ref{prop:Casorati} 
\[
\big(\mu(\ga_1)-\mu(\ga_2)\big) \langle \phi, \overline{\psi} \rangle_{k;n} =D(\phi, \psi)(z_+q^{n-1})- D(\phi, \psi)(z_-q^{k-1}). \qedhere
\]
\end{proof}
We are going to apply the previous lemma to the functions $\varphi_\ga$, $\varphi_\ga^\dagger$ and $\Phi_\ga^\pm$, which are functions in $V_{\mu(\ga)}$ by Propositions \ref{prop:basis Vmu} and \ref{prop:Phi+-}. The following lemma will be useful.
\begin{lem} \label{lem:asympt DPhi}
For $k \rightarrow \infty$,
\[
D(\Phi_{\ga_1}^\pm,\Phi_{\ga_2}^\pm)(z_\pm q^{-k}) = (\ga_1-\ga_2)  K_{z_\pm}(\ga_1 \ga_2)^{k-1}\big(1+\mathcal O(q^k)\big).
\]
\end{lem}
\begin{proof}
This follows from the definition of the Casorati determinant \eqref{def:Casorati}, and from the asymptotic behavior of $\Phi_\ga^\pm(x)$ and $u(x)/x$ for large $|x|$, see \eqref{eq:asymp Phi} and Lemma \ref{lem:asymp w&u}. See also the proof of Lemma \ref{lem:det Phi}.
\end{proof}
As a consequence we obtain the following orthogonality relation.
\begin{lem} \label{lem:orth varphi Phi}
Let $\ga \in \T$ and $\ga' \in \Ga$, then
\[
\langle \varphi_\ga, \Phi^+_{\ga'} \rangle_{\mathcal L^2} =0, \qquad \langle \varphi_\ga^\dagger, \Phi^+_{\ga'} \rangle_{\mathcal L^2} =0.
\]
\end{lem}
\begin{proof}
From Lemmas \ref{lem:inprod phi psi}, \ref{lem:asympt DPhi} and the $c$-function expansions from Proposition \ref{prop:c function exp} we find
\[
\langle \varphi_\ga, \Phi^+_{\ga'} \rangle_{\mathcal L^2} = \lim_{k,n \rightarrow -\infty}\big( c_{z_+}(\ga) \langle \Phi^+_\ga, \Phi^+_{\ga'} \rangle_{k;n} +c_{z_+}(1/\ga)\langle \Phi^+_{1/\ga}, \Phi^+_{\ga'} \rangle_{k;n} \big) = 0,
\]
since $|\ga'|<1$. In the same way it follows that $\langle \varphi_\ga^\dagger, \Phi^+_{\ga'} \rangle_{\mathcal L^2} =0$.
\end{proof}

We are now ready to show that $\mathcal G_p$ is a partial right inverse of the map $\mathcal F$.
\begin{prop} \label{prop:Gp}
The map $\mathcal G_p : \mathcal H_p \rightarrow F(\R_q)$ satisfies
\[
\langle \mathcal G_p g_1, \mathcal G_pg_2 \rangle_{\mathcal L_2} = \langle g_1, g_2 \rangle_{\mathcal H_p}, \qquad g_1,g_2 \in \mathcal H_p.
\] 
Moreover, for $g \in \mathcal H_p$ we have $\mathcal F (\mathcal G_{p}g) = \mathbf 0 +g$ in $\mathcal H$, where $\mathbf 0$ denotes the zero function in $\mathcal H_c$.
\end{prop}
\begin{proof}
Let $g,h$ be finitely supported functions in $\mathcal H_p$. Then we find from Corollary \ref{cor:orth rel},
\[
\begin{split}
\langle \mathcal G_pg, \mathcal G_ph \rangle_{\mathcal L^2} &= \int_{\R_q} \Big( \sum_{\ga \in \Ga} g(\ga) \Phi_\ga^+(x) N(\ga) \Big)  \overline{\Big( \sum_{\ga' \in \Ga} h(\ga') \Phi_{\ga'}^+(x) N(\ga') \Big)}w(x) d_qx\\
&= \sum_{\ga,\ga' \in \Ga} \langle \Phi_\ga^+,  \Phi_{\ga'}^+ \rangle_{\mathcal L^2} g(\ga) \overline{h(\ga')} N(\ga) N(\ga') \\
&= \sum_{\ga \in \Ga} g(\ga) \overline{h(\ga)} N(\ga) = \langle g, h \rangle_{\mathcal H_p}. 
\end{split}
\]
In order to prove the identity $\mathcal F(\mathcal G_p g)=\mathbf 0 + g$, we split this identity into three different cases: 
\[
\mathcal F_p (\mathcal G_p g)=g, \qquad \mathcal F_c (\mathcal G_p g)=0 \quad \text{and} \quad \mathcal F_c^\dagger (\mathcal G_p g)=0.
\]
The identity $\mathcal F_p (\mathcal G_p g)= g$ is proved in a similar way as in the proof of Proposition \ref{prop:GF=id}, and the other two identities follow from Lemma \ref{lem:orth varphi Phi}. Since the set of finitely supported functions is dense in $\mathcal H_p$, the proposition follows.
\end{proof}
Next we are going to show that $\mathcal G_c$ is also a partial right inverse of $\mathcal F$. For this we apply a classical method used by G\"otze \cite{Go65} and by Braaksma and Meulenbeld \cite{BM67} for the Jacobi function transform.

We define for $\ga \in \T$, 
\begin{equation} \label{eq:def u1 u2}
\begin{split}
u_1(\ga)=u_1(\ga;a,b,c,d;z_-,z_+|q)&=K_{z_+}\,c_{z_+}(\ga)c_{z_+}(1/\ga)-K_{z_-}\,c_{z_-}(\ga)c_{z_-}(1/\ga),\\
u_2(\ga)=u_2(\ga;a,b,c,d;z_-,z_+|q)& = K_{z_+}\,c_{z_+}(\ga)c_{z_+}^\dagger(1/\ga)-K_{z_-}\,c_{z_-}(\ga)c_{z_-}^\dagger (1/\ga).
\end{split}
\end{equation}
Explicitly, using the expressions for $c_z$ and $K_z$, we have
\[
\begin{split}
u_1(\ga) &=   \frac{(1-q) (s\ga^{\pm 1}, cq\ga^{\pm 1}/as, dq\ga^{\pm 1}/as;q)_\infty }{ (cq/a,cq/a,dq/a,dq/a,\ga^{\pm 2};q)_\infty \te(bz_+, bz_-, cz_+, cz_-, dz_+, dz_-)} \\
& \times \Big( z_+\te(az_+, bz_-, cz_-, dz_-, bsz_+\ga^{\pm 1})-z_- \te(az_-, bz_+, cz_+, dz_+, bsz_-\ga^{\pm 1}) \Big).
\end{split}
\]
For $u_2$ we have
\[
\begin{split}
u_2(\ga) =& \frac{(1-q)q\ga}{as}\frac{ (s\ga^{\pm 1}, cq/as\ga, cq\ga/bs, dq/as\ga, dq/bs;q)_\infty }{ (cq/a, cq/b, dq/a, dq/b, \ga^{\pm2};q)_\infty } \\
&\times \Bigg( \frac{\te(bsz_-\ga,asz_-/q\ga)}{\te(cz_-,dz_-)} - \frac{\te(bsz_+\ga,asz_+/q\ga)}{\te(cz_+,dz_+)}  \Bigg).
\end{split}
\]
Using the $\te$-product identity \eqref{eq:theta prod} with
\[
x=z_-e^{i(\kappa+\de)/2} \sqrt{|cd|}, \quad y= z_+ e^{i(\kappa+\de)/2}\sqrt{|cd|}, \quad v = \frac{bs \ga e^{-i(\kappa+\de)/2}}{\sqrt{|cd|}}, \quad w= e^{i(\kappa-\de)/2}\sqrt{\left| \frac{c}{d} \right| },
\]
where $c=|c|e^{i\kappa}$ and $d= |d|e^{i\de}$, we obtain
\[
u_2(\ga) =  z_+(1-q) \frac{ (s\ga^{\pm1},  cq\ga^{\pm1}/as, cq\ga^{\pm1}/bs , dq\ga^{\pm 1}/as, dq\ga^{\pm1}/bs;q)_\infty \te(z_-/z_+, cdz_-z_+)}{(cq/a, cq/b, dq/a, dq/b, \ga^{\pm2};q)_\infty  \te(cz_-, cz_+, dz_+, dz_-) }.
\]
Observe that $u_2=u_2^\dagger$ and $u_2(\ga)=u_2(1/\ga)$, and that $u_2$ is real-valued on $\T\setminus\{-1,1\}$.

Let $C_0(\T)$ be the set of functions defined by
\[
C_0(\T) = \Big\{ g:\T \rightarrow \C \ | \ g \text{ is continuous,}\ g(-1)=g(1)=0,\ g(\ga)=g(1/\ga) \Big\}.
\]  
\begin{prop} \label{prop:lim kn}
Let $g \in C_0(\T)$ and let $\ga'\in \T\setminus\{-1,1\}$, then
\[
\begin{split}
\lim_{k,n \rightarrow -\infty} \frac{1}{4\pi i} \int_\T g(\ga) \langle \varphi_\ga, \overline{\varphi_{\ga'}} \rangle_{k;n}\frac{d\ga}{\ga } &= g(\ga')u_1(\ga'), \\
\lim_{k,n \rightarrow -\infty} \frac{1}{4\pi i} \int_\T g(\ga) \langle \varphi_\ga^\dagger, \overline{\varphi_{\ga'}} \rangle_{k;n}\frac{d\ga}{\ga } &=g(\ga')u_2(\ga') 
\end{split}
\]
\end{prop}
\begin{proof}
We prove the first identity in the proposition, the second identity is proved in the same way. 
Let us fix a $g \in C_0(\T)$, and let us define
\[
I_{z}^m(\te') = \frac{1}{2\pi} \int_0^{\pi} g(e^{i\te}) \frac{ D(\varphi_{e^{i\te}}, \varphi_{e^{i\te'}})(zq^{-m})} { 2\cos(\te) - 2 \cos(\te') } d\te.
\]
From Lemma \ref{lem:inprod phi psi} we find
\[
\frac{1}{4\pi i}\int_\T g(\ga) \langle \varphi_\ga, \overline{\varphi_{\ga'}} \rangle_{k;n}\frac{d\ga}{\ga} = I_{z_+}^{1-n}(\te') - I_{z_-}^{1-k}(\te'),
\]
where $\ga'= e^{i\te'}$ with $\te' \in (0,\pi)$. We see that we need to investigate the limit of $I_z^m(\te')$ when $m \rightarrow \infty$. Using the $c$-function expansion from Proposition \ref{prop:c function exp} and Lemma \ref{lem:asympt DPhi} we obtain, for large $m$,
\[
I_z^m(\te') = \frac{K_z}{2\pi} \sum_{\xi,\eta \in \{-1,1\}} \int_0^\pi g(e^{i\te}) \Big(\psi_z^m(\te,\te';\xi,\eta) +\mathcal O(q^m) \Big) d\te,
\]
where
\[
\psi_z^m(\te,\te';\xi,\eta) = \frac{ (e^{i \xi \te}- e^{i \eta \te'}) e^{i(m-1)(\xi\te+\eta \te')}c_z(\xi \te) c_z(\eta \te') }{2 \cos(\te) - 2 \cos(\te')}.
\]
Since $c_z(\ga)$ is continuous on $\T\setminus\{-1,1\}$ the functions $\psi_z^m$, considered as functions of $\te$, are continuous on $(0,\pi)\setminus\{\te'\}$. We see immediately that $\psi_z^m(\te,\te';1,1)$ and $\psi_z^m(\te,\te';-1,-1)$ have a removable singularity at $\te=\te'$. Now by the Riemann-Lebesgue Lemma the terms with these two functions vanish in the limit, and this leaves us with
\[
\lim_{m \rightarrow \infty} I_z^m(\te') = \lim_{m \rightarrow \infty} \frac{K_z}{2\pi}\int_0^\pi g(e^{i\te}) \big(\psi_z^m(\te,\te';1,-1)+ \psi_z^m(\te,\te';-1,1)\big) d\te.
\] 
Here we applied dominated convergence to get rid of the $\mathcal O(q^{m})$-terms.
Using the identity $\cos(\al) - \cos(\be) = 2 \sin(\frac{\al+\be}{2}) \sin(\frac{ \al-\be}{2})$ we find
\begin{align}
\psi_z^m(\te,\te';1,-1)+ \psi_z^m(\te,\te';-1,1) & \nonumber\\
=\frac{ 1}{4 \sin\left(\frac{\te+\te'}{2} \right) \sin\left( \frac{ \te - \te'}{2}\right) } &\Big(c_z(e^{i\te}) c_z(e^{-i\te'}) e^{i(m-1)(\te-\te')} (e^{i\te}-e^{-i\te'}) \nonumber \\
\qquad &+ c_z(e^{-i\te}) c_z(e^{i\te'}) e^{i(m-1)(\te'-\te)} (e^{-i\te}-e^{i\te'}) \Big) \nonumber \\
=\frac{ 1}{4 \sin\left(\frac{\te+\te'}{2} \right) \sin\left( \frac{ \te - \te'}{2}\right) } &
\Big( [c_z(e^{-i\te})c_z(e^{i\te'})- c_z(e^{i\te})c_z(e^{-i\te'})] e^{i(m-1)(\te'-\te)}(e^{-i\te}-e^{i\te'}) \nonumber\\
&+ c_z(e^{i\te}) c_z(e^{-i\te'}) \psi^m(\te,\te') \Big), \label{eq:psi}
\end{align}
where 
\[
\begin{split}
\psi^m(\te,\te') &= e^{i(m-1)(\te'-\te)}(e^{-i\te}-e^{i\te'}) + e^{i(m-1)(\te-\te')}(e^{i\te}-e^{-i\te'})\\
&= 2\cos\big(m\te -(m-1) \te'\big) - 2\cos\big((m-1)\te - m\te'\big).
\end{split}
\]
The first term in \eqref{eq:psi} has a removable singularity, so by the Riemann-Lebesgue Lemma this term also vanishes in the limit, and now we have
\[
\begin{split}
\lim_{m \rightarrow \infty} I_z^m(\te')& = \lim_{m \rightarrow \infty} \frac{K_z}{2\pi}\int_0^\pi g(e^{i\te})c_z(e^{i\te}) c_z(e^{-i\te'})  \frac{\psi^m(\te,\te') }{4 \sin\left(\frac{\te+\te'}{2} \right) \sin\left( \frac{ \te - \te'}{2}\right) }d\te\\
&= \lim_{m \rightarrow \infty} \frac{K_z}{2\pi}\int_0^\pi g(e^{i\te})c_z(e^{i\te}) c_z(e^{-i\te'}) D_m(\te;\te') d\te,
\end{split}
\]
where $D_m(\te;\te')$ is the Dirichlet kernel
\[
D_m(\te;\te') = \frac{ \sin\left( (m-\hf)(\te-\te') \right) }{ \sin\left(\hf(\te-\te') \right) }.
\]
From the well-known properties of the Dirichlet-kernel we obtain
\[
\lim_{m \rightarrow \infty} I_z^m(\te') = K_z\, g(e^{i\te'}) c_z(e^{i\te'})c_z(e^{-i\te'}),
\]
and from this the result follows.
\end{proof}

\begin{prop} \label{prop:Fg1g2}
Let $g_1,g_2 \in C_0(\T)$ and let $\ga' \in \T\setminus\{-1,1\}$, then
\[
\int_{\R_q} \Bigg[ \frac{1}{4\pi i} \int_\T \vect{\varphi_\ga(x)}{\varphi_\ga^\dagger(x)}^T \vect{g_1(\ga)}{g_2(\ga)} \frac{ d\ga}{\ga}\Bigg] \vect{\varphi_{\ga'}(x)}{\varphi_{\ga'}^\dagger(x)} w(x) d_qx= \mathbf{u}(\ga')\vect{g_1(\ga')}{g_2(\ga')},
\]
where $\mathbf u$ is the matrix-valued function on $\T\setminus\{-1,1\}$ defined by
\[
\ga \mapsto \mathbf{u}(\ga) =
\begin{pmatrix}
u_2(\ga) & u_1(\ga) \\
u_1^\dagger(\ga) & u_2(\ga)
\end{pmatrix}.
\]
\end{prop}
\begin{proof}
Let $g_1,g_2 \in C_0(\T)$ and $\ga,\ga' \in \T\setminus\{-1,1\}$. From Proposition \ref{prop:lim kn} we find
\[
\begin{split}
g_2(\ga')u_1(\ga') &= \lim_{k,n \rightarrow -\infty} \frac{1}{4\pi i} \int_\T g_2(\ga) \Big(\int_{z_-q^{k}}^{z_+q^{n}} \varphi_\ga(x) \varphi_{\ga'}(x)w(x) d_qx\Big) \frac{d\ga}{\ga} \\
&= \int_{\R_q} \varphi_{\ga'}(x)\Big( \frac{1}{4\pi i} \int_\T
g_2(\ga) \varphi_\ga(x) \frac{d\ga}{\ga} \Big) w(x) d_qx.
\end{split}
\]
To justify the interchanging of the order of integration, we note that it follows from the explicit expressions for $\varphi_\ga(x)$ and $w(x)$ that
\[
\varphi_\ga(x) \varphi_{\ga'}(x)w(x) = 1 + \mathcal O(q^m), \qquad x=z_{\pm}q^m \rightarrow 0,
\]
so that the sums
\[
\sum_{m=n}^\infty \varphi_\ga(z_\pm q^m) \varphi_{\ga'}(z_\pm q^m) q^m w(z_\pm q^m), 
\]
both converge uniformly on $\T\setminus\{-1,1\}$. In the same way we find
\[
\begin{split}
g_2(\ga')u_2(\ga') = \int_{\R_q} \varphi_{\ga'}^\dagger(x)\Big( \frac{1}{4\pi i} \int_\T
g_2(\ga) \varphi_\ga(x) \frac{d\ga}{\ga} \Big) w(x) d_qx,\\
g_1(\ga')u_1^\dagger(\ga') = \int_{\R_q} \varphi_{\ga'}^\dagger(x)\Big( \frac{1}{4\pi i} \int_\T
g_1(\ga) \varphi_\ga^\dagger(x) \frac{d\ga}{\ga} \Big) w(x) d_qx, \\
g_1(\ga')u_2(\ga') = \int_{\R_q} \varphi_{\ga'}(x)\Big( \frac{1}{4\pi i} \int_\T
g_1(\ga) \varphi_\ga^\dagger(x) \frac{d\ga}{\ga} \Big) w(x) d_qx.
\end{split}
\]
Now the proposition follows.
\end{proof}
The matrix-valued function $\mathbf u$ has the following useful property, which is proved in the appendix.
\begin{lem} \label{lem:u^-1 = v}
For $\ga \in \T \setminus \{-1,1\}$,
\[
\mathbf{u}(\ga)^{-1} = \mathbf{v}(\ga).
\]
\end{lem}
We define
\[
C_0(\T;\C^2) = \Bigg\{ g=\vect{g_1}{g_2}\ \Big|\ g_1,g_2 \in C_0(\T) \Bigg\} \subset \mathcal H_c.
\]
We are now in a position to show that $\mathcal G_c$ is a partial right inverse of $\mathcal F$.
\begin{prop} \label{prop:Gc}
The map $\mathcal G_c:\mathcal H_c \rightarrow F(\R_q)$ satisfies
\[
\langle \mathcal G_cg_1, \mathcal G_cg_2 \rangle_{\mathcal L^2} = \langle g_1,g_2 \rangle_{\mathcal H_c}, \qquad g_1,g_2 \in \mathcal H_c.
\] 
Moreover, for $g \in \mathcal H_c$ we have $\mathcal F(\mathcal G_c g)=g+\mathbf 0$ in $\mathcal H$, where $\mathbf 0$ denotes the zero function in $\mathcal H_p$.
\end{prop}
\begin{proof}
Let $\ga' \in \T\setminus\{-1,1\}$, let $g^{(1)},g^{(2)} \in C_0(\T)$ and define $g= \vect{g^{(1)}}{g^{(2)}}$. Since $v_1$, $v_1^\dagger$ and $v_2$ are continuous on $\T$, both components of the $\C^2$-valued function
\[
\ga \mapsto
\mathbf{v}(\ga)
\vect{g^{(1)}(\ga)}{g^{(2)}(\ga)}
\]
are in $C_0(\T)$. Now by Proposition \ref{prop:Fg1g2} and Lemma \ref{lem:u^-1 = v} we have 
\[
\vect{\big(\mathcal F_c(\mathcal G_c g)\big)(\ga')}{\big(\mathcal F_c^\dagger(\mathcal G_c g)\big)(\ga')} = \mathbf u(\ga') \mathbf v(\ga')g(\ga') = g(\ga'). 
\]
Moreover, for $\ga' \in \Ga$,
\[
\begin{split}
\big(\mathcal F_p ( \mathcal G_c g)\big)(\ga') &= \lim_{k,n \rightarrow -\infty} \langle \mathcal G_c g, \Phi_{\ga'}^+ \rangle_{k;n} \\
&= \lim_{k,n \rightarrow -\infty} \frac{1}{4\pi i} \int_\T \big( \langle \varphi_\ga^\dagger, \Phi_{\ga'}^+ \rangle_{k;n} \ \langle \varphi_\ga, \Phi_{\ga'}^+ \rangle_{k;n} \big) \mathbf v(\ga) g(\ga) \frac{d\ga}{\ga} \\
&=0,
\end{split}
\]
by dominated convergence and Lemma \ref{lem:orth varphi Phi}. This shows that $\mathcal F(\mathcal G_c g) = g+\mathbf 0$ in $\mathcal H$. 

Let $g_1,g_2 \in C_0(\T^2)$, then it follows from Proposition \ref{prop:F iso} that
\[
\langle g_1,g_2 \rangle_{\mathcal H_c} =\langle g_1+\mathbf 0, g_2+\mathbf 0 \rangle_{\mathcal H}= \langle \mathcal F(\mathcal G_c  g_1), \mathcal F(\mathcal G_c g_2) \rangle_{\mathcal H} = \langle \mathcal G_c g_1, \mathcal G_c g_2 \rangle_{\mathcal L^2}.
\]
Since the set $C_0(\T;\C^2)$ is dense in $\mathcal H_c$, the proposition follows.
\end{proof}

Collecting the results of this subsection we come to the main theorem.
\begin{thm} \label{thm:F}
For $(a,b,c,d) \in P$, the map $\mathcal F:\mathcal L^2 \rightarrow \mathcal H$ is an isometric isomorphism with inverse $\mathcal G$.
\end{thm}
\begin{proof}
Let $(a,b,c,d) \in P_{\mathrm{gen}}$. Combining Propositions \ref{prop:Gp} and \ref{prop:Gc} gives $\mathcal F \mathcal G= id_{\mathcal H}$. Together with Proposition \ref{prop:GF=id} this leads to the theorem.  By continuity in the parameters, the result holds for all $(a,b,c,d) \in P$.
\end{proof}
\begin{cor} \label{cor:orth basis}
The set $\{ \Psi(x,\cdot) \ | \ x \in \R_q\}$ forms an orthogonal basis for $\mathcal H$ with squared norm $\|\Psi(x,\cdot) \|^2_{\mathcal H} = w(x)^{-1}$.
\end{cor}
\begin{proof}
This follows from Lemma \ref{lem:Psi in H} and Theorem \ref{thm:F}, since the functions $f_y$ defined in Lemma \ref{lem:Psi in H} form an orthogonal basis for $\mathcal L^2$ with squared norm $w(y)^{-1}$.
\end{proof}
\begin{rem}
The Hilbert space $\mathcal H$ and the inverse $\mathcal G$ of the vector-valued big $q$-Jacobi function transform depend essentially on five parameters, namely
$az_-,bz_-,cz_-,dz_-$ and $z_+/z_-$.
\end{rem}
\subsection{An equivalent integral transform}
For $f \in \mathcal D_{\mathrm{fin}}$ and $\ga \in \Ga$ we have $(\mathcal Ff)(\ga) = \langle f, \Phi_\ga^+ \rangle_{\mathcal L^2}$. Since the function $\Phi_\ga^+$ can be expressed in terms of big $q$-Jacobi functions by Proposition \ref{prop:Phi+-}, we can define an integral transform with only the big $q$-Jacobi functions $\varphi_\ga$ and $\varphi_\ga^\dagger$ as a kernel, which is equivalent to $\mathcal F$. This new integral transform can of course also be extended to an isometric isomorphism. We only state the result here, and we leave the details to the reader.\\
 
For $f \in \mathcal D_{\mathrm{fin}}$ we define an integral transform $\mathcal J$ that is closely related to the vector-valued big $q$-Jacobi function transform $\mathcal F$ by
\[
(\mathcal J f)(\ga) = \int_{\R_q} f(x) \vect{\varphi_\ga(x)}{\varphi_\ga^\dagger(x)} w(x) d_q x, \qquad \ga \in \T\cup \Ga.
\]
For $(a,b,c,d) \in P$ we define an inner product on the vector space of $\C^2$-valued functions by
\[
\langle f,g \rangle_{\mathcal M} = \frac{ 1}{4\pi i } \int_\T g(\ga)^T 
\mathbf v(\ga)
f(\ga) \frac{ d\ga}{\ga} + \sum_{\ga \in \Ga} g(\ga)^T \mathbf v_p(\ga) f(\ga).
\]
Here $\mathbf v_p(\ga)$ is the matrix-valued function on $\Ga$ given by
\[
\ga \mapsto \mathbf v_p(\ga)=
\begin{pmatrix}
v_{3,p}(\ga) & v_{4,p}(\ga) \\
v_{1,p}(\ga) & v_{2,p}(\ga) 
\end{pmatrix},
\]
where the matrix coefficients $v_{i,p}(\ga)=v_{i,p}(\ga;a,b,c,d;z_-,z_+|q)$, $i=1,\ldots,4$, are defined as follows:\\
For $\ga \in \Ga^{\mathrm{inf}} \cup \Ga^{\mathrm{fin}}_{q/s}$,  $v_{4,p}(\ga) = v_{1,p}^\dagger(\ga)$, $v_{3,p}(\ga) = v_{2,p}(\ga)$, and
\begin{align*}
v_{1,p}(\ga)&= d_{z_+}^2(\ga) N(\ga), \\
v_{2,p}(\ga) &= d_{z_+}(\ga)d_{z_+}^\dagger(\ga)N(\ga).
\end{align*}
For $\ga \in \Ga^{\mathrm{fin}}_{dq/as} \cup \Ga^{\mathrm{fin}}_{s}$, $v_{2,p}(\ga)=v_{4,p}(\ga)=0$, and
 
\[
\begin{split}
v_{1,p}(\ga)& =
\begin{cases}
0, & \text{if } a=\overline{b},\\
\dfrac{ N(\ga)}{\big(c_{z_+}(\ga)\big)^2 }, & \text{if } a \neq \overline{b},
\end{cases} \\
v_{2,p}(\ga) & =
\begin{cases}
\dfrac{ N(\ga)}{c_{z_+}(\ga) c_{z_+}^\dagger(\ga) }, & \text{if } a=\overline{b},\\
0, & \text{if } a \neq \overline{b}.
\end{cases}
\end{split}
\]
Recall here that $\Ga^{\mathrm{fin}}_{dq/as}$ is only non-empty  if $a\neq\overline{b}$. Now denote by $\mathcal M= \mathcal M(a,b,c,d;z_-,z_+|q)$ the closure of the set
\[
\mathrm{span} \left\{ \ga \rightarrow \vect{ \varphi_\ga(x) }{ \varphi_\ga^\dagger(x) } \ \Big| \ x \in \R_q \right\}
\]
with respect to the norm $\| \cdot \|_{\mathcal M}$. Note that a function $g \in \mathcal M$ satisfies $rg=g$.

Let $\Theta:\mathcal M \rightarrow \mathcal H$ be the operator defined by
\[
(\Theta g)(\ga) = 
\begin{cases}
g(\ga), & \ga \in \T,\\
\big(d_{z_+}(\ga) \ d_{z_+}^\dagger(\ga) \big) g(\ga),& \ga \in \Ga^{\mathrm{inf}} \cup \Ga^{\mathrm{fin}}_{q/s},\\
\big( c_{z_+}(\ga)^{-1}\ 0 \big) g(\ga), & \ga \in \Ga^{\mathrm{fin}}_s\cup \Ga^{\mathrm{fin}}_{dq/as},
\end{cases}
\]
then 
\[
\langle \Theta g_1, \Theta g_2 \rangle_{\mathcal H} = \langle g_1, g_2 \rangle_{\mathcal M}, 
\]
for functions  $g_i \in \mathcal M$, $i=1,2$. In particular, we have
\[
\left(\Theta \vect{\varphi_{\cdot}(x) }{\varphi_{\cdot}^\dagger(x)}\right)(\ga) = \Psi(x,\ga), \qquad x \in \R_q, \quad \ga \in \T \cup \Ga,
\]
so $\Theta:\mathcal M \rightarrow \mathcal H$ is an isomorphism. Also, $\mathcal Ff = (\Theta \circ \mathcal J)f$ for $f \in \mathcal D_{\mathrm{fin}}$.

\begin{thm} 
The map $\mathcal J: \mathcal D_{\mathrm{fin}} \rightarrow \mathcal M$ extends uniquely to an isometric isomorphism $\mathcal J_{\mathrm{ext}}:\mathcal L^2 \rightarrow \mathcal M$. Moreover, $\mathcal I = \mathcal G \circ \Theta : \mathcal M \rightarrow \mathcal L^2$ is the inverse of $\mathcal J_{\mathrm{ext}}$.
\end{thm}

\begin{rem}
(i) Let $f \in \mathcal L^2$ be a function for which $\mathcal Ff$ can be written as an integral transform, i.e.,
\[
(\mathcal Ff)(\ga) = 
\int_{\R_q} f(x) \Psi(x,\ga) w(x) d_q x, \qquad \ga \in \T \cup \Ga.
\] 
Then $\mathcal J_{\mathrm{ext}}f$ can in general \emph{not} be written as the integral
\[
\int_{\R_q} f(x) \vect{\varphi_\ga(x) }{ \varphi_\ga^\dagger(x)} w(x) d_q x,
\]
when $f \not\in \mathcal D_{\mathrm{fin}}$, since the integrals in the components of this vector might be divergent for $\ga \in \Ga$.

(ii) The inverse $\mathcal I$ of $\mathcal J_{\mathrm{ext}}$ can be given explicitly by
\[
(\mathcal Ig)(x) = \Bigg \langle g, \vect{\varphi_\cdot(x)}{\varphi_{\cdot}^\dagger(x)} \Bigg \rangle_{\mathcal M}, \qquad x \in \R_q, 
\]
for all $g \in \mathcal M$ for which the above inner product exists.

(iii) Like $\mathcal F$, the map $\mathcal J_{\mathrm{ext}}$ diagonalizes $L$;
\[
(\mathcal J_{\mathrm{ext}} \circ L \circ \mathcal J_{\mathrm{ext}}^{-1})g = M g,
\]
for functions $g \in \mathcal M$ such that $M g \in \mathcal M$.
\end{rem}

\appendix
\section{}
In this appendix we prove Lemmas \ref{lem:I(x,y)} and \ref{lem:u^-1 = v}.

\subsection{Proof of Lemma \ref{lem:I(x,y)}} \label{appA1}
We prove the following statement:\\
\emph{For $x,y \in \R_q$ and  $\ga, \ga^{-1} \in \mathcal S_{\mathrm{reg}}\setminus \mathcal V$ we have
\[
\begin{split}
&\frac{ \Phi_{1/\ga}^-(x) \Phi_{1/\ga}^+(y)}{v(1/\ga)} - \frac{ \Phi_{\ga}^-(x) \Phi_{\ga}^+(y)}{v(\ga)} = \\ 
&\quad\frac{1}{\ga-1/\ga}\Big[ v_1(\ga) \varphi_\ga(x) \varphi_\ga(y) + v_2(\ga) \big(\varphi_\ga(x) \varphi_\ga^\dagger(y) + \varphi_\ga^\dagger(x) \varphi_\ga(y) \big) + v_1^\dagger(\ga) \varphi_\ga^\dagger(x)\varphi_\ga^\dagger(y)\Big],
\end{split}
\]
where
\begin{align*}
v_1(\ga) &= \frac{ (cq/a, dq/a;q)_\infty^2 \te(bz_+,bz_-) }{(1-q)ab z_-^2z_+^2 \te(z_-/z_+,z_+/z_-,a/b,b/a)}\\
& \times \frac{(\ga^{\pm2};q)_\infty}{(s\ga^{\pm1}, cq\ga^{\pm1}/as, dq\ga^{\pm1}/as;q)_\infty \te(s\ga^{\pm1},absz_- z_+\ga^{\pm 1})}\\
& \times \Big( z_-\te(az_+, cz_+, dz_+, bz_-, asz_-\ga^{\pm1}) -  z_+\te(az_-, cz_-, dz_-, bz_+, asz_+\ga^{\pm1})  \Big),\\
v_2(\ga) &= \frac{(cq/a,dq/a,cq/b,dq/b;q)_\infty \te(az_+, az_-, bz_+, bz_-, cdz_- z_+) }{abz_-^2z_+(1-q)\te(z_+/z_-,a/b,b/a)} \\
& \times \frac{ (\ga^{\pm2};q)_\infty }{ (s\ga^{\pm1};q)_\infty \te(s\ga^{\pm1}, absz_- z_+\ga^{\pm1}) }.
\end{align*} 
}
\begin{proof}
Let $\ga, \ga^{-1} \in \mathcal S_{\mathrm{reg}}\setminus \mathcal V$. Note that $\mathcal S_{\mathrm{pol}} \subset \mathcal V $, hence $\ga, \ga^{-1} \not\in \mathcal S_{\mathrm{pol}}$. We define
\[
I_\ga(x,y) = \frac{ \Phi_{1/\ga}^-(x) \Phi_{1/\ga}^+(y)}{v(1/\ga)} - \frac{ \Phi_{\ga}^-(x) \Phi_{\ga}^+(y)}{v(\ga)}.
\]
Using $\varphi_\ga= \varphi_{1/\ga}$ and Proposition \ref{prop:Phi+-} we see that
\[
I_\ga(x,y) = v_1'(\ga) \varphi_\ga(x)\varphi_\ga(y) +v_2'(\ga) \varphi_\ga(x)\varphi_\ga^\dagger(y)+ v_3'(\ga) \varphi_\ga^\dagger(x)\varphi_\ga(y)+ v_4'(\ga)\varphi_\ga^\dagger(x)\varphi_\ga^\dagger(y),
\]
where
\begin{align*}
v_1'(\ga) &= \frac{ d_{z_-}(1/\ga)d_{z_+}(1/\ga) }{ v(1/\ga) } - \frac{ d_{z_-}(\ga)d_{z_+}(\ga) }{ v(\ga) },\\
v_2'(\ga) &= \frac{ d_{z_-}(1/\ga)d_{z_+}^\dagger(1/\ga) }{ v(1/\ga)} - \frac{ d_{z_-}(\ga)d_{z_+}^\dagger(\ga) }{ v(\ga) },\\
v_3'(\ga) &= \frac{ d_{z_-}^\dagger(1/\ga)d_{z_+}(1/\ga) }{ v(1/\ga) } - \frac{ d_{z_-}^\dagger(\ga)d_{z_+}(\ga) }{ v(\ga) },\\
v_4'(\ga) &= \frac{ d_{z_-}^\dagger(1/\ga)d_{z_+}^\dagger(1/\ga) }{ v(1/\ga) } - \frac{ d_{z_-}^\dagger(\ga)d_{z_+}^\dagger(\ga) }{ v(\ga) }.
\end{align*}
Since $v(\ga)=v^\dagger(\ga)$, it is immediately clear that $v_1'^\dagger(\ga) = v_4'(\ga)$ and $v_2'^\dagger(\ga)=v_3'(\ga)$. 

Using the explicit expressions for $d_z(\ga)$ and $v(\ga)$, see Proposition \ref{prop:Phi+-} and Theorem \ref{thm:properties Phi}, we find
\[
\begin{split}
v_2'(\ga) =& \frac{bs (cq/a, dq/a, cq/b, dq/b;q)_\infty \te(az_+,bz_-) }{q(1-q)(s\ga, s/\ga;q)_\infty \te(z_-/z_+,a/b,b/a)}\\
& \times \Bigg(\frac{ \te(q^2/asz_-\ga, q/bsz_+\ga) }{ \te(s\ga,q^2/absz_- z_+\ga) } - \frac{ \te(q^2\ga/asz_-, q\ga/bsz_+) }{ \te(s/\ga,q^2\ga/absz_- z_+) }  \Bigg)
\end{split}
\]
From this we find the expression for $v_2(\ga) = (\ga-1/\ga) v_2'(\ga)$ given in the lemma after using the $\te$-product identity \eqref{eq:theta prod} with 
\begin{align*}
x&=\frac{ iqe^{-i\al/2} }{s}\sqrt{\left|\frac{q }{az_-}\right|}, & y &=  \frac{i qe^{-i\al/2}}{bsz_+} \sqrt{ \left| \frac{q}{az_-} \right|}, \\
 v&=  i \ga e^{-i\al/2}  \sqrt{\left|\frac{q}{az_-} \right|}, & w &= \frac{ie^{-i\al/2}}{\ga} \sqrt{\left|\frac{q} {az_-} \right|},
\end{align*}
where $a=|a|e^{i\al}$.

Next we compute $v_1(\ga) = (\ga-1/\ga) v_1'(\ga)$;
\[
\begin{split}
v_1'(\ga) =&\frac{ (cq/a, dq/a;q)_\infty^2 \te(bz_+,bz_-) }{\ga z_+(1-q)(s\ga,s/\ga;q)_\infty \te(z_-/z_+) \te(a/b)^2} \\
 \times& \Bigg(\frac{\ga^2 (cq\ga/bs, dq\ga/bs;q)_\infty \te(q^2\ga/ as z_+, q^2\ga/ asz_-) }{ (cq\ga/as, dq\ga/as;q)_\infty \te(s/ \ga, qs\ga/cdz_- z_+)} \\
& -\frac{ (qc/bs\ga , dq/bs\ga;q)_\infty \te(q^2/\ga as z_+, q^2/\ga asz_-) }{ (cq/as \ga, dq/ad\ga;q)_\infty \te(s \ga, qs/cdz_- z_+\ga)}\Bigg).
\end{split}
\]
Since $cq/as=bs/c$, the expression between large brackets can be written as
\begin{equation} \label{eq:expression1}
\begin{split}
&\frac{-\ga}{(cq\ga^{\pm1}/as, dq\ga^{\pm1}/as;q)_\infty \te(s\ga^{\pm1}, cdz_- z_+\ga^{\pm1}/s)} \\
\times& \Big(\ga^{-1} \te(\ga asz_-/q, \ga asz_+/q, dq/sb\ga, cq/sb\ga, s/\ga, cdz_- z_+/s\ga)\\
&\quad -\ga\te(asz_-/q\ga, asz_+/q\ga, dq\ga/bs, cq\ga/bs, s\ga, cdz_- z_+ \ga/s) \Big).
\end{split}
\end{equation}
We use the $\te$-product identity \cite[Exer.5.22]{GR04}
\[
\begin{split}
\frac{1}{y}\te(tx/p, ux/p, vx/p, wx/p, y/p, y/r, r/p)-\frac{1}{x}\te(ty/p, uy/p, vy/p, wy/p, x/p, x/r, r/p) = \\
\frac{1}{y} \te(tr/p, ur/p, vr/p, wr/p, x/p, y/p, y/x) -\frac{x}{pr} \te(t,u,v,w, r/x, r/y,y/x),
\end{split}
\]
with parameters
\[
\begin{split}
&p=\frac{q}{asz_+},\quad r=\frac{q}{asz_-}, \quad t=\frac{q}{az_+}, \quad  u=\frac{q}{cz_+}, \\
&v=\frac{q}{dz_+},\quad  w= bz_-,\quad   x=\frac{1}{\ga}, \quad y=\ga,
\end{split}
\]
then \eqref{eq:expression1} becomes
\[
\begin{split}
&\frac{-\te(\ga^2)}{(cq\ga^{\pm1}/as, dq\ga^{\pm1}/as;q)_\infty \te(s\ga^{\pm1}, cdz_- z_+\ga^{\pm1}/s,z_+/z_-)}\\
& \times \Big(
\te(q/az_-, q/cz_-, q/dz_-, bz_+, asz_+\ga^{\pm1}/q) \\
& \quad- \frac{a^2s^2z_- z_+}{q^2} \te(q/az_+, q/cz_+, q/dz_+, bz_-, q\ga^{\pm1}/asz_-)  \Big).
\end{split}
\]
The expression given in the lemma is obtained from this after using the identity $-x\te(qx)=\te(x)$ several times.
\end{proof}

\subsection{Proof of Lemma \ref{lem:u^-1 = v}}
We show that 
\[
\mathbf u(\ga)^{-1}=
\mathbf v(\ga), \qquad \ga \in \T\setminus\{-1,1\},
\]
with
\[
\mathbf u(\ga)=\begin{pmatrix}
u_2(\ga) & u_1(\ga) \\
u_1^\dagger(\ga) & u_2(\ga)
\end{pmatrix}, \qquad 
\mathbf v(\ga) =
\begin{pmatrix}
v_2(\ga) & v_1^\dagger(\ga) \\
v_1(\ga) & v_2(\ga)
\end{pmatrix}.
\]
\begin{proof}
By a direct verification, using the explicit expressions for $v_1$ and $v_2$ from Lemma \ref{lem:I(x,y)} (see also Appendix \ref{appA1}), and for $u_1$ and $u_2$, one sees that
\[
v_2(\ga) =\frac{u_2(\ga)}{\de(\ga)}, \qquad v_1^\dagger(\ga) =-\frac{u_1(\ga)}{\de(\ga)},
\]
where $\de(\ga)$ is the function given by
\[
\begin{split}
\de(\ga) &=  \frac{(1-q)^2 z_-^2 z_+^2 ab\, \te(z_-/z_+, z_+/z_-, a/b, b/a) }{(cq/a, cq/b, dq/a, dq/b;q)_\infty^2 \te(az_-, az_+, bz_-, bz_+, cz_-, cz_+, dz_-, dz_+,)} \\
& \times \frac{(s\ga^{\pm1}, s\ga^{\pm1}, cq\ga^{\pm1}/as,  cq\ga^{\pm1}/bs, dq\ga^{\pm1}/as, cq\ga^{\pm1}/bs)_\infty \te(s\ga^{\pm1}, absz_- z_+ \ga^{\pm1})}{(\ga^{\pm2};q)_\infty^2 }
\end{split}
\]
It remains to show that $\de(\ga)$ is the determinant of the matrix $\mathbf u(\ga)$.

Using the definition \eqref{eq:def u1 u2} of the functions $u_1$ and $u_2$, and $u_2=u_2^\dagger$, the determinant of $\mathbf u(\ga)$ becomes
\[
\begin{split}
\det\big(\mathbf u(\ga)\big) =& u_2^\dagger(\ga) u_2(\ga) - u_1^\dagger(\ga) u_1(\ga) \\
=& K_{z_-}K_{z_+} \Big( c_{z_+}(\ga) c_{z_+}(1/\ga) c_{z_-}^\dagger(\ga) c_{z_-}^\dagger(1/\ga)
-c_{z_+}(\ga) c_{z_+}^\dagger(1/\ga) c_{z_-}^\dagger(\ga) c_{z_-}(1/\ga) \\
&\phantom{K_{z_-}K_{z_+}} +c_{z_+}^\dagger(\ga) c_{z_+}^\dagger(1/\ga) c_{z_-}(\ga) c_{z_-}(1/\ga)
-c_{z_+}^\dagger(\ga) c_{z_+}(1/\ga) c_{z_-}(\ga) c_{z_-}^\dagger(1/\ga) \Big)\\
=& K_{z_-}K_{z_+} \Big(c_{z_+}(\ga) c_{z_-}^\dagger(\ga) - c_{z_+}^\dagger(\ga) c_{z_-}(\ga) \Big) \\
& \times \Big(c_{z_+}(1/\ga) c_{z_-}^\dagger(1/\ga) - c_{z_+}^\dagger(1/\ga) c_{z_-}(1/\ga) \Big).
\end{split}
\]
Explicitly, we have
\[
c_{z_+}(\ga) c_{z_-}^\dagger(\ga) - c_{z_+}^\dagger(\ga) c_{z_-}(\ga) = \frac{ (s/\ga,s/\ga, cq/as\ga, cq/bs\ga, dq/as\ga, dq/bs\ga;q)_\infty F(\ga)}{ (cq/a, cq/b, dq/a, dq/b,1/\ga^2, 1/\ga^2;q)_\infty^2 \te(az_-, az_+, bz_-, bz_+)} ,
\]
where
\[
\begin{split}
F(\ga) =& \te(asz_-\ga, bsz_+\ga, bz_-, az_+) - \te(asz_+\ga, bsz_-\ga, bz_+, az_-) \\
=& bz_+ \te(z_-/z_+, a/b, s\ga, absz_-z_+\ga).
\end{split}
\]
See the proof of Theorem \ref{thm:properties Phi} for the last equality. By inspection it follows that indeed $\de(\ga) = \det\big(\mathbf u(\ga)\big)$.
\end{proof}


\begin{thebibliography}{99}
\bibitem{AA85} G.E. Andrews, R. Askey, Classical orthogonal polynomials. in: \textit{ Polyn\^omes orthogonaux et applications.}, Lecture Notes in Math., 1171, Springer, Berlin, 1985, 36-62.

\bibitem{AW85} R. Askey, J. Wilson, \textit{Some basic hypergeometric
orthogonal polynomials that generalize Jacobi polynomials}, Mem. Amer. Math. Soc. \textbf{54} (1985), no. 319.

\bibitem{BM67} B.L.J. Braaksma, B. Meulenbeld, \textit{Integral transforms with generalized Legendre functions as kernels}, Compositio Math. \textbf{18} (1967), 235-287.

\bibitem{DS63} N. Dunford, J.T. Schwartz, \textit{Linear Operators Part II}, Interscience, New York, 1963. 

\bibitem{GR04} G. Gasper, M. Rahman, \textit{Basic Hypergeometric Series}, 2nd ed., Cambridge University Press, Cambridge, 2004.

\bibitem{Go65} F. G\"otze, \textit{Verallgemeinerung einer Integraltransformation von Mehler-Fock durch den von Kuipers und Meulenbeld eingef\"urten Kern $P_k^{m,n}(z)$}, Indag. Math. \textbf{27} (1965), 396-404. 

\bibitem{GKR05} W. Groenevelt, E. Koelink, H. Rosengren, \textit{Continuous Hahn functions as Clebsch-Gordan coefficients}, in: Theory and applications of special functions, M.E.H. Ismail, E. Koelink (Eds.), Dev. Math., \textbf{13}, Springer, New York, 2005, 221-284.

\bibitem{GIM96} D.P. Gupta, M.E.H. Ismail, D.R. Masson, \textit{Contiguous relations, basic hypergeometric functions, and orthogonal polynomials. III. Associated continuous dual $q$-Hahn polynomials},  J. Comput. Appl. Math. \textbf{68} (1996), 115-149.

\bibitem{K95} T. Kakehi, \textit{Eigenfunction expansion associated with the Casimir operator on the quantum group ${\rm SU}_q(1,1)$}, Duke Math. J. \textbf{80} (1995), 535-573.

\bibitem{KMU95} T. Kakehi, T. Masuda, K. Ueno, \textit{Spectral analysis of a $q$-difference operator which arises from the quantum ${\rm SU}(1,1)$ group},  J. Operator Theory \textbf{33} (1995), 159-196.

\bibitem{KS98} R. Koekoek, R.F. Swarttouw, \textit{The Askey-scheme of hypergeometric orthogonal polynomials and its q-analogue}, Report
98-17, Technical University Delft, Delft, 1998.

\bibitem{KS01} E. Koelink, J.V. Stokman, \textit{The Askey-Wilson function transform}, Internat. Math. Res. Notices (2001), 1203-1227.

\bibitem{KS03} E. Koelink, J.V. Stokman, \textit{The big $q$-Jacobi function transform},  Constr. Approx. \textbf{19} (2003), 191-235.

\bibitem{Koo84} T.H. Koornwinder, \textit{Jacobi functions and analysis on noncompact semisimple Lie groups}, in: Special Functions: Group Theoretical Aspects and Applications, R.A. Askey, T.H. Koornwinder, W. Schempp (Eds.), D. Reidel Publ. Comp., Dordrecht, 1984, 1-85.

\bibitem{Ma75} R.P. Martin,\textit{On the decomposition of tensor products of principal series representations for real-rank one semisimple groups},  Trans. Amer. Math. Soc. \textbf{201} (1975), 177-211.

\bibitem{Ne05} Yu.A. Neretin, \textit{Some Continuous Analogs of the Expansion in Jacobi Polynomials and Vector-Valued Orthogonal Bases},  Funct. Anal. Appl. \textbf{39} (2005), no. 2, 106-119.

\bibitem{Pu61} L. Puk\'anszky, \textit{On the Kronecker products of irreducible representations of the $2\times 2$ real unimodular group. I,} Trans. Amer. Math. Soc. \textbf{100} (1961), 116-152.

\bibitem{We10} H. Weyl, \textit{\"Uber gew\"ohnliche lineare Differentialgleichungen mit singulare Stellen und ihre Eigenfunktionen, (2. Note)}, Nachr. K\"onig. Gesell. Wiss. G\"ottingen. Math.-Phys. Kl., (1910), 442-467.
\end{thebibliography}
\end{document}